%% file: main.tex
\documentclass[reqno]{amsart}
\usepackage[utf8]{inputenc}

\usepackage{fullpage}

\usepackage{amsmath}
\usepackage{amssymb}
\usepackage[foot]{amsaddr}
\usepackage{bm}
\usepackage{mathdots}

\usepackage{graphicx}
\usepackage{subcaption}

\usepackage{algorithm}
\usepackage{algpseudocode}

\usepackage[colorlinks, citecolor=blue]{hyperref}

\newcommand{\fref}[1]{Fig.\,\ref{#1}}

\newcommand{\eref}[1]{Eq.\,(\ref{#1})}

\newcommand{\sref}[1]{Sec.\!~\ref{#1}}

\newcommand{\cref}[1]{Ref.\,\cite{#1}}

\newcommand{\pref}[1]{Alg.\,\ref{#1}}

\let\Vec\relax
\newcommand{\Vec}[1]{\bm{#1}}
\newcommand{\Mat}[1]{\mathbf{#1}}
\newcommand{\RR}{\mathbb{R}}
\newcommand{\CCGrid}[1]{\mathsf{CCGrid}(#1)}
\newcommand{\tr}{\operatorname{tr}}

\DeclareMathOperator{\Uniform}{Uniform}

\newcommand{\ScalarParamNds}[1]{\lambda^\text{nds}_{#1}}                          
\newcommand{\ParamNds}[1]{\Vec{\lambda}^\text{nds}_{#1}}                          
\newcommand{\StateNdsExact}[2][t]{\Vec{\phi}^\text{nds}_{#2}(#1)}                 
\newcommand{\StateNdsApprx}[2][t]{\Vec{\hat{\phi}}^\text{nds}_{#2}(#1)}           
\newcommand{\DynamicsNdsExact}[2][\Vec{\phi}]{\Vec{r}^\text{nds}_{#2}(#1)}        
\newcommand{\DynCoeffNdsExact}[1]{\Mat{A}^\text{nds}_{#1}}                        

\newcommand{\ScalarParamSmp}[1]{\lambda^\text{smp}_{#1}}                          
\newcommand{\ParamSmp}[1]{\Vec{\lambda}^\text{smp}_{#1}}                          
\newcommand{\StateSmpExact}[2][t]{\Vec{\phi}^\text{smp}_{#2}(#1)}                 
\newcommand{\DynamicsSmpExact}[2][\Vec{\phi}]{\Vec{r}^\text{smp}_{#2}(#1)}        

\newcommand{\StateSmpStdSC}[2][t]{\Vec{\hat{\phi}}^\text{s-sc}_{#2}(#1)}          

\newcommand{\StateSmpDynSC}[2][t]{\Vec{\hat{\phi}}^\text{d-sc}_{#2}(#1)}          
\newcommand{\DynamicsSmpDynSC}[2][\Vec{\phi}]{\Vec{\hat{r}}^\text{d-sc}_{#2}(#1)} 
\newcommand{\DynCoeffSmpDynSC}[1]{\Mat{\hat{A}}^\text{d-sc}_{#1}}                 

\title{Accurate Data-Driven Surrogates of Dynamical Systems for Forward Propagation of Uncertainty}
\author{Saibal De$^{1,\dagger}$, Reese E.\ Jones$^{2,\dagger}$, and Hemanth Kolla$^{3,\dagger}$}
\address{ORCID ID: $^1$0000-0003-4691-189X, $^2$0000-0002-2332-6279, $^3$0000-0003-4969-5870}
\address{$^\dagger$Sandia National Laboratories, Livermore, CA 94551}

\begin{document}

\begin{abstract}
Stochastic collocation (SC) is a well-known non-intrusive method of constructing surrogate models for uncertainty quantification.
In dynamical systems, SC is especially suited for full-field uncertainty propagation that characterizes the distributions of the high-dimensional primary solution fields of a model with stochastic input parameters.
However, due to the highly nonlinear nature of the parameter-to-solution map in even the simplest dynamical systems, the constructed SC surrogates are often inaccurate.
This work presents an alternative approach, where we apply the SC approximation over the dynamics of the model, rather than the solution.
By combining the data-driven sparse identification of nonlinear dynamics (SINDy) framework with SC, we construct dynamics surrogates and integrate them through time to construct the surrogate solutions.
We demonstrate that the SC-over-dynamics framework leads to smaller errors, both in terms of the approximated system trajectories as well as the model state distributions, when compared against full-field SC applied to the solutions directly.
We present numerical evidence of this improvement using three test problems: a chaotic ordinary differential equation, and two partial differential equations from solid mechanics.
\end{abstract}

\maketitle
\input{sections/introduction}
\input{sections/background}
\input{sections/dynamics_sc}
\input{sections/results}
\input{sections/summary}

\section*{Acknowledgements}
This work is supported by the Laboratory Directed Research and Development (LDRD) program at Sandia National Laboratories.
Sandia National Laboratories is a multi-mission laboratory managed and operated by National Technology and Engineering Solutions of Sandia, LLC (NTESS), a wholly owned subsidiary of Honeywell International Inc., for the U.S. Department of Energy’s National Nuclear Security Administration (DOE/NNSA) under contract DE-NA0003525. This written work is authored by an employee of NTESS. The employee, not NTESS, owns the right, title and interest in and to the written work and is responsible for its contents. Any subjective views or opinions that might be expressed in the written work do not necessarily represent the views of the U.S. Government. The publisher acknowledges that the U.S. Government retains a non-exclusive, paid-up, irrevocable, world-wide license to publish or reproduce the published form of this written work or allow others to do so, for U.S. Government purposes. The DOE will provide public access to results of federally sponsored research in accordance with the DOE Public Access Plan.

\bibliographystyle{unsrt}
\bibliography{references}

\end{document}

%% file: sections/introduction.tex

\section{Introduction}

Uncertainty quantification (UQ) is an essential component of robust model-based predictions \cite{le2010spectral,anjos2012handbook,smith2013uncertainty,ghanem2017handbook}, as variations in model inputs (parameters, initial or boundary conditions etc.) can drastically change the response of a system.
This is especially relevant for high-consequence applications, such as engineered systems in extreme environments and autonomous mobility, where slight deviations can lead to catastrophic failures.
In complex physical systems and related  mathematical models, these points of failures can originate from multiple, sometimes unexpected, sources.
Characterizing the behavior of a few pre-determined quantities of interest (QoIs) is therefore not sufficient in these scenarios.
To mitigate this limitation, a previous work \cite{jones2021minimally} from our group has advocated treating the primary solution fields (e.g.\ displacements and velocities of a solid body in a structural mechanics simulation) as the QoIs, and characterizing the uncertainty of these full-field solutions.

Given a characterization of parameter uncertainty, forward uncertainty propagation involves pushing these uncertainties forward to uncertainties of QoIs. 
In scientific computing applications, this is generally achieved via Monte Carlo sampling.
A large ensemble of independent forward simulations, initialized with different values of uncertain model inputs, is required to offset the slow convergence rate of the Monte Carlo QoI estimates.
Due to the high computational cost of each forward solve for typical engineering applications, UQ studies typically rely on cheap surrogate models of the input-to-QoIs map.

For full-field uncertainty propagation, the target QoIs can be significantly high-dimensional; this limits the type of surrogates we can use.
Regression based surrogates, such as neural networks \cite{tripathy2018deep, raissi2019physics, cai2021physics} and Gaussian processes \cite{rasmussen2005gaussian, giovanis2020data}, can be difficult and computationally expensive to train.
Furthermore neural networks are well-known to be data hungry, while Gaussian processes have issues with large inputs.
Projection-based surrogates, such as polynomial chaos expansions \cite{debusschere2017intrusive}, involve computing inner product integrals over the stochastic input parameter space.
Interpolation based surrogates, such as stochastic collocation \cite{xiu2017stochastic}, are non-intrusive and relatively inexpensive, especially for high-dimensional uncertain inputs.

In this work, we focus on stochastic collocation (SC) surrogates; their construction requires only the ability to compute QoIs at given system inputs.
SC evaluates the full model at suitably placed collocation points in the parameter space, then computes the surrogates using interpolation from these exact values.
This black-box treatment of the computational model, coupled with the relatively low cost of interpolation even with high-dimensional QoIs, makes SC highly attractive for full-field UQ.
Several review articles \cite{xiu2007efficient, eldred2009recent, jakeman2009stochastic, xiu2009fast, xiu2017stochastic} provide comprehensive descriptions of the method, including the interpolation and regression formulations of SC, and construction of sparse collocation grids \cite{smolyak1963quadrature, barthelmann2000high, garcke2013sparse} for high-dimensional parameter spaces.

Despite the significant advantages of stochastic collocation surrogates in terms of computational cost and ease of implementation, their fidelity to the full model can be low in high-dimensional input spaces.
When the response of a system is highly non-linear in terms of the model inputs, the method requires a relatively dense collocation grid to accurately approximate the QoIs, which increases the number of full model evaluations and total cost.
This is especially relevant when approximating the full solution state of a system of differential equations with uncertain parameters.
In many cases the parameter-to-state map can be highly complex even when the driving term has a simple form.
A remedy, proposed previously by us in \cref{jones2021minimally}, is based on the hypothesis that the dynamics of a system of differential equations has a simpler dependence on the uncertain parameters, compared to the dependence of solution state on the same uncertain parameters.
Accordingly, directly approximating the dynamics, rather than the solution state, within the stochastic collocation framework will drastically improve the accuracy of the surrogates.
However, in \cref{jones2021minimally} stochastic collocation was employed to express the approximation of dynamics purely as a function of uncertain parameters, which limited the accuracy improvements.
Specifically, for certain problem classes where the SC interpolation and numerical integration operators commute, this offered no accuracy benefits over standard SC applied to solution state. 

The main contribution of the present article is the development of an approach that pushes the dynamics approximation paradigm considerably beyond that of \cref{jones2021minimally}.
The central limitation of \cref{jones2021minimally} was not capturing the implicit dependence of the dynamics on the solution state.
We address it here by \emph{expressing the approximate dynamics explicitly as function of both solution state and uncertain parameters}.
With SC as our choice for efficient approximation along the stochastic dimensions, we adopt the sparse identification of nonlinear dynamics (SINDy) framework \cite{brunton2016discovering} to capture the dynamics dependence on solution state.
With judicious choice of SINDy basis, we construct an approximation of the dynamics as a separable function of the solution state and uncertain parameters, relying on the expressivity of the SINDy framework for accuracy and the efficiency of SC for feasibility.
The overall formulation is applicable to a broad class ordinary differential equations (ODE) systems, e.g. predator-prey models, chemical reacting systems, control problems, epidemiology; but our primary application is the discretized dynamics of deformable solids governed by large-scale partial differential equations (PDEs). 

The rest of the paper is structured as follows.
In \sref{sec:background}, we formalize the notation for forward uncertainty quantification in the context of dynamical systems and briefly review the stochastic collocation methodology.
We also summarize the SINDy framework for dynamics discovery \cite{brunton2016discovering}, a tool that we use in \sref{sec:dynamics_sc} to develop our proposed method.
\sref{sec:results} contains results from three demonstration problems: (a) the chaotic Lorenz ODE system, (b) a one-dimensional solid bar impact problem, and (c) a quasi three-dimensional multi-material notch impact problem;
the first problem pertains to an ODE system while the second and third problems treat discretized PDE systems from solid mechanics. 
In each case, we compare the accuracy our dynamics surrogates to that of the standard stochastic collocation surrogates over the state space.
Finally, \sref{sec:summary} summarizes our findings and lists several avenues for future work.


%% file: sections/background.tex

\section{Background} \label{sec:background}

Our formulation of forward uncertainty propagation is motivated by an approach that solves probability density function (PDF) equations associated with the differential equation system under uncertainty.
We first present a brief overview of the PDF evolution equation and outline a sampling based approach that makes its solution computationally tractable.
We then review how using SC surrogates accelerates this sampling task.
Next, we adapt this approach for PDEs discretized with finite elements, specifically dynamic solid mechanics problems.
We conclude this section with a brief review of SINDy, which allows us to learn the unknown driving terms of a dynamical systems in a data-driven fashion.
To reiterate, our proposal is to perform accurate UQ by approximating the dependence of the dynamics, rather than the system state, with respect to the model parameters within a combined SC and SINDy framework, which we present in the \sref{sec:dynamics_sc}.

\subsection{Forward Uncertainty Propagation} \label{sec:forwarduq}

A large class of science and engineering applications involve evolving a system, from some initial state, $\Vec{\phi}_0$, according to prescribed dynamics:
\begin{equation}
  \label{eq:ode_system}
  \dot{\Vec{\phi}} = \Vec{r}(\Vec{\phi}, \Vec{\lambda}) + \Vec{r}_\text{ext}(t, \Vec{\lambda}, \Vec{\phi}), \quad \Vec{\phi}(t = 0, \Vec{\lambda}) = \Vec{\phi}_0.
\end{equation}
Here $\Vec{\lambda}$ are parameters characterizing the model dynamics, $\Vec{\phi}(t, \Vec{\lambda}) \in \RR^S$ is the state of the system at time $t$, $\dot{\Vec{\phi}} \equiv \partial_t\Vec{\phi}$ is the time derivative of the state, $\Vec{r}(\Vec{\phi}, \Vec{\lambda}) \in \RR^S$ is the time-independent internal model dynamics, and $\Vec{r}_\text{ext}(t, \Vec{\phi}, \Vec{\lambda}) \in \RR^S$ encapsulates any time-dependent external forcing.
Computation of the external force is generally much simpler and cheaper than evaluating the internal forces. 
In this work, we ignore the external forcing term for the purpose of method development ($\Vec{r}_\text{ext} \equiv \Vec{0}$).
Also note that the state vector $\Vec{\phi}$ is often obtained by discretizing a PDE (see \sref{sec:ode_from_pde}), so that it encapsulates variations over the spatial dimensions.

Very few model parameters are exactly known; they are either measured with physical instruments with limited precision or obtained from inference studies that assign uncertainties to the predicted values.
Mathematically, we model this by treating $\Vec{\lambda}$ as a random variable with a distribution determined by experimental data.
We assume that the parameters do not evolve in time (i.e.\ $\dot{\Vec{\lambda}} = \Vec{0}$), and that they vary in a piecewise constant manner through the solid domain (e.g.\ when a solid object is made out of multiple materials).
Then the parameters are completely characterized by a finite dimensional vector ($\Vec{\lambda} \in \RR^P$), and the stochasticity is fully captured by a joint PDF $p_{\Vec{\lambda}}(\Vec{\lambda})$.
Note that we are using the same symbols for a random variable as well as its instantiation; we will use this convention throughout the rest of this article.

Stochasticity in the model parameters, $\Vec{\lambda}$, implies that the model state $\Vec{\phi}$ is also a random variable, with an associated time-varying PDF $p_{\Vec{\phi}}(\Vec{\phi}; t)$.
The same holds true for any QoIs $\Vec{q} = \Vec{q}(\Vec{\phi})$ derived from the model state $\Vec{\phi}$.
Forward uncertainty propagation studies characterize the uncertainties in the QoIs given the form of the parameter PDF $p_{\Vec{\lambda}}$ and the model dynamics $\Vec{r}$.
For full-field UQ, rather than tracking a few pre-determined QoIs over time, we directly learn the state PDF $p_{\Vec{\phi}}$; from this distribution, we can then obtain the distributions of any derived QoI.
This eliminates the need to repeat computations when investigating the behavior of a new QoI.

\subsection{The PDF Evolution Equation}\label{sec:PDEeq}

Pope \cite{pope1985lagrangian} derives the exact evolution equation for the joint state-parameter PDF associated with \eref{eq:ode_system} as:
\begin{equation}
  \label{eq:pdf_evolution}
  \partial_t p_{\Vec{\phi} \Vec{\lambda}}(\Vec{\phi}, \Vec{\lambda}; t) = -
  \nabla_{\Vec{\phi}} \cdot \{p_{\Vec{\phi} \Vec{\lambda}}(\Vec{\phi},
  \Vec{\lambda}; t) \Vec{r}(\Vec{\phi}, \Vec{\lambda})\}.
\end{equation}
The initial state-parameter distribution at time $t = 0$ is given by: $p_{\Vec{\phi} \Vec{\lambda}}(\Vec{\phi}, \Vec{\lambda}; 0) = \delta(\Vec{\phi} - \Vec{\phi}_0) p_{\Vec{\lambda}}(\Vec{\lambda})$, where $\delta$ represents the Dirac-$\delta$ distribution describing a point mass density.
The state distribution is then obtained by marginalizing the joint distribution over the model parameters:
\begin{equation}
  \label{eq:pdf_state}
  p_{\Vec{\phi}}(\Vec{\phi}; t) = \int p_{\Vec{\phi} \Vec{\lambda}}(\Vec{\phi},
  \Vec{\lambda}; t) \mathrm{d}\Vec{\lambda}.
\end{equation}
In obtaining the state distribution numerically, we encounter two challenges: first, the PDE in \eref{eq:pdf_evolution} is very high dimensional, and traditional discretization methods run into the ``curse of dimensionality'', requiring exponentially more computational resources as the dimensionality increases.
Second, the integral in \eref{eq:pdf_state} is also high-dimensional for systems with a large number of parameters, so that numerical estimates may be expensive and inaccurate.

A statistically equivalent approach \cite{pope1985lagrangian} is to draw samples $\{\ParamSmp{i} : 1 \leq i \leq I\} \sim p_{\Vec{\lambda}}$ from the parameter distribution, solve \eref{eq:ode_system} to generate corresponding system trajectories $\StateSmpExact{i} \equiv \Vec{\phi}(t, \ParamSmp{i})$, and use these sample trajectories to approximate the state distribution (e.g.\ using kernel density estimation).
This way, the cost of solving the high-dimensional PDF evolution equation is avoided, and replaced with sampling and forward simulations of the original ODE system.

These Monte-Carlo type sampling-based approaches have become the \emph{de facto} standard for forward UQ of computational models.
They exhibit the familiar $O(1 / \sqrt{I})$ convergence rate for estimating quantities of interest.
Given this slow rate, and the typical high computational cost of running forward real-scale simulations, most practitioners rely on computationally cheaper surrogate models to accelerate UQ studies.

\subsection{Stochastic Collocation} \label{sec:sc_state}

Stochastic collocation (SC) \cite{eldred2009recent} is one of the simplest surrogate modeling strategies for approximating a quantity of interest\footnote{%
Note the slight abuse notations here; the QoI $\Vec{q}$ depends on state $\Vec{\phi}$, which in turn depends on parameters $\Vec{\lambda}$.
Written in this form, we are suppressing the intermediate state variable.}
$\Vec{q}(\Vec{\lambda})$ with stochastic parameters $\Vec{\lambda}$.
Relying on a few exact evaluations $\{\Vec{q}(\ParamNds{j}) : 1 \leq j \leq J\}$ of the QoI at collocation nodes $\ParamNds{j}$, it constructs approximations at the sample parameters $\ParamSmp{i}$ through interpolation: 
\begin{equation}
  \label{eq:sc}
  \Vec{q}(\ParamSmp{i}) \approx \sum_{j = 1}^{J} \gamma_{ij} \, \Vec{q}(\ParamNds{j}), \quad 1 \leq i \leq I.
\end{equation}

In systems with one uncertain parameter $\lambda \in \RR$, the collocation points are traditionally chosen using standard quadrature rules (e.g.\ the Gauss-Legendre grid) to enable computations of moments of the QoI distributions.
In our experiments, we use the nested Clenshaw-Curtis grids \cite{clenshaw1960method}.
To simplify the exposition and without loss of generality, assume that the range of the bounded parameter $\lambda$ is restricted to the interval $[-1, 1]$.
In this standard interval, the level-$\ell$ Clenshaw-Curtis grid is defined as:
\begin{equation}
  \label{eq:clenshaw_curtis_1d}
  \CCGrid{\ell} = \bigg\{\ScalarParamNds{j} = \cos\frac{(j - 1) \pi}{2^\ell} \; \bigg\vert \; 1 \leq j \leq 2^\ell + 1\bigg\}, \quad \ell \in \{0, 1, 2, \ldots\}.
\end{equation}
Using these nodes, we approximate a QoI as:
\begin{equation}
  \Vec{q}(\lambda) \approx \sum_{j = 1}^{2^\ell + 1} \gamma_j(\lambda) \Vec{q}(\ScalarParamNds{j}), \quad \gamma_j(\ScalarParamNds{j'}) =
  \begin{cases}
    1 & \text{if} \quad j = j', \\
    0 & \text{if} \quad j \neq j',
  \end{cases}
\end{equation}
and for a given sample parameter value $\ScalarParamSmp{i}$, the interpolation coefficients are computed as $\gamma_{ij} = \gamma_j(\ScalarParamSmp{i})$.
A popular choice for the basis functions $\gamma_j$ are the Lagrange polynomials:
\begin{equation}
  \gamma_j(\lambda) = \prod_{\substack{j' = 1\\j' \neq j}}^{2^\ell + 1} \frac{\lambda - \ScalarParamNds{j'}}{\ScalarParamNds{j} - \ScalarParamNds{j'}};
\end{equation}
they allow exact representation of any polynomials in the parameter $\lambda$ with maximum degree $2^\ell$.
Since the grids are nested, with $\CCGrid{\ell} \subseteq \CCGrid{\ell + 1}$ for all $\ell \geq 0$, we can easily construct adaptively refined interpolation schemes; this is a major advantage over traditional Gaussian quadratures.

In higher dimensions, with $\Vec{\lambda} \in [-1, 1]^P$, a natural choice of collocation points is constructed by taking tensor products of the one-dimensional grids:
\begin{equation}
  \begin{split}
    \CCGrid{P, \Vec{\ell}}
    &= \CCGrid{\ell_1} \otimes \cdots \otimes \CCGrid{\ell_P} \\
    &= \{\Vec{\lambda}^\text{nds}_{j_1, \cdots, j_P} = (\ScalarParamNds{j_1}, \ldots, \ScalarParamNds{j_P}) : 1 \leq j_p \leq 2^{\ell_p} + 1, 1 \leq p \leq P\};
  \end{split}
\end{equation}
the corresponding basis functions $\gamma_{j_1, \ldots, j_P}$ are obtained by multiplying the one-dimensional functions:
\begin{equation}
  \gamma_{j_1, \ldots, j_P}(\Vec{\lambda}) = \gamma_{j_1}(\lambda_1) \cdots \gamma_{j_P}(\lambda_P), \quad \Vec{\lambda} = (\lambda_1, \ldots, \lambda_P).
\end{equation}

\input{figures/sparse_grid/demo}

In the tensor grid, the number of collocation points grows exponentially with increasing number of stochastic parameters, and exactly evaluating the QoI at all the nodes become prohibitively costly.
Sparse grids \cite{smolyak1963quadrature} can significantly reduce the number of nodes by strategically eliminating nodes associated with high-order polynomials, as demonstrated in \fref{fig:sparse_grid_demo}; hence, sparse grids are regularly employed for approximating high-dimensional functions across a wide range of applications \cite{barthelmann2000high, garcke2013sparse}.

As noted in an earlier work \cite{jones2021minimally}, SC with sparse grids is particularly effective for constructing surrogate trajectories in dynamical systems.
To see this, denote the exactly simulated trajectories corresponding to collocation nodes as $\StateNdsExact{j} = \Vec{\phi}(t, \ParamNds{j})$; then \eref{eq:sc} specializes to:
\begin{equation}
  \label{eq:sc_state}
  \StateSmpExact{i} \approx \StateSmpStdSC{i} = \sum_{j = 1}^{J} \gamma_{ij} \, \StateNdsExact{j}, \quad 1 \leq i \leq I.
\end{equation}
Here the `hat' in $\StateSmpStdSC{i}$ indicates that this quantity is an approximation, and the superscript `s-sc' indicates that we used stochastic collocation over the state space (State SC) to construct it.
Note that the coefficients $\gamma_{ij}$ do not change with time; they can be be pre-computed after the sample and collocation parameters have been chosen, and reused throughout the rest of the simulation.
This feature allows us to construct an ``online'' approximation scheme, where we only access the exactly simulated system states from the current time step.
This can drastically reduce the memory footprint of the forward propagation method.
For this reason, State SC is an attractive framework for full-field UQ of dynamical systems.

\subsection{Finite Element Discretization of PDEs} \label{sec:ode_from_pde}

Many ODE systems in engineering applications are constructed by discretizing PDEs using traditional finite difference, volume, or element schemes.
To provide a basis for our methodological developments, we briefly review the finite element method (FEM) in the context of continuum solid mechanics.

The primary governing equation of solid mechanics is the linear momentum balance:
\begin{equation} \label{eq:solid}
  \rho_0 \ddot{\Vec{u}} = \nabla_{\Vec{X}} \cdot \Mat{P}.
\end{equation}
Here $\Vec{X}$ is the reference location of a point in the solid body $\Omega \subseteq \RR^D$, $\Vec{x} = \Vec{x}(t, \Vec{X})$ is the location of the reference point at time $t$, $\Vec{u}(t, \Vec{X}) = \Vec{x} - \Vec{X}$ is the displacement, $\Mat{P}$ is the first Piola-Kirchhoff stress, $\rho_0$ is the density of the material, and $D$ is the dimensionality of the enclosing space.
In formulating the FEM discretization, we consider the weak form of the governing equations, obtained by integrating the residual of \eref{eq:solid} against a test function $\Vec{\omega}(\Vec{X})$:
\begin{equation}
  \int_\Omega \Vec{\omega} \cdot \left(\nabla_{\Vec{X}} \cdot \Mat{P} -
  \rho_0 \ddot{\Vec{u}}\right) \mathrm{d}\Vec{X}= 0.
\end{equation}
We select the test function from an appropriate Sobolev space with basis functions $\omega_n(\Vec{X})$, indexed by the label $n \in \{1, \ldots, N\}$ of the corresponding finite element node $\Vec{X}_n$.
Integrating the residual along the span of any of these basis functions by parts, we obtain the Galerkin projection:
\begin{equation}
  \int_\Omega \left(- \nabla \omega_n \cdot \Mat{P} - \rho_0 \omega_n
  \ddot{\Vec{u}}\right) \mathrm{d}\Vec{X} + \int_{\partial\Omega} \omega_n
  \Mat{P} \cdot \mathrm{d}\Vec{S} = 0.
\end{equation}
Using the same basis to express the displacement, $\Vec{u}(t, \Vec{X}) = \sum_{n'} \Vec{u}_{n'}(t) \omega_{n'}(\Vec{X})$, we obtain:
\begin{equation}
  \sum_j \underbrace{\left(\int_\Omega \rho_0 \omega_n \omega_{n'} \mathrm{d}\Vec{X}\right)}_{M_{n n'}} \ddot{\Vec{u}}_{n'}
  =
  \underbrace{-\int_\Omega \nabla \omega_n \cdot \Mat{P} \mathrm{d}\Vec{X} + \int_{\partial \Omega} \omega_n \Mat{P} \cdot \mathrm{d}\Vec{S}}_{\Vec{f}_n}.
\end{equation}
and thereby a matrix system of coupled ODEs:
\begin{equation}
  \ddot{\Vec{\phi}} = \Vec{r},
  \quad
  \Vec{\phi} = \begin{bmatrix} \vdots \\ \Vec{u}_n \\ \vdots \end{bmatrix},
  \quad
  \Vec{r} = \Mat{M}^{-1} \Vec{f},
  \quad
  \Vec{f} = \begin{bmatrix} \vdots \\ \Vec{f}_n \\ \vdots \end{bmatrix},
  \quad
  \Mat{M} = \begin{bmatrix} \ddots & \vdots & \iddots \\ \cdots & M_{n n'} \Mat{I}_D & \cdots \\ \iddots & \vdots & \ddots
  \end{bmatrix},
\end{equation}
where $\Mat{I}_D$ is the $D \times D$ identity matrix.
The dimensionality $S = D N$ of the discretized system state can easily number in thousands to hundreds of thousands for even a moderately sized problem.
Evaluating the dynamics $\Vec{r}$ becomes increasingly costly with complex constitutive models for the stress $\Mat{P}$, higher-order FEM basis functions, and/or finer FEM mesh sizes.

\subsection{Sparse Identification of Nonlinear Systems}

The costly evaluation of the dynamics $\Vec{r}$ is the main bottleneck in constructing an accurate SC surrogate.
We alleviate this computational burden by treating the exact dynamics as ``unknown'' and learning a cheaper surrogate dynamics from system trajectories in a data-driven manner.
The sparse identification of nonlinear systems (SINDy) framework, introduced in Ref. \cite{brunton2016discovering}, is a systematic way of learning unknown dynamics of ODE systems from trajectory data.
It has also been extended to learn and analyze implicit ODEs \cite{mangan2016inferring}, PDEs \cite{schaeffer2017learning} and ODEs with stochastic parameters \cite{corbetta2020application}, across various fields such as nonlinear control \cite{brunton2016sparse}, biological networks \cite{mangan2019model}, and fluid mechanics \cite{loiseau2018sparse}.

Let us briefly summarize the key components of the original SINDy: we express the unknown dynamics as a sparse linear combination of $B$ basis functions\footnote{%
We are suppressing the dependence of the dynamics on model parameters in this section to simplify notation.}:
\begin{equation}
  r_s(\Vec{\phi}) \approx \sum_{b = 1}^B A_{s b} \psi_b(\Vec{\phi}), \quad 1 \leq s \leq S,
\end{equation}
where $A_{s b}$ is the coefficient describing the contribution of the basis $\psi_b$ to the $s$-th component of the dynamics.
We rewrite the dynamics approximation matrix form as:
\begin{equation}
  \label{eq:sindy_factorization}
  \Vec{r}(\Vec{\phi}) \approx \Mat{A} \, \Vec{\psi}(\Vec{\phi}), \quad
  \Mat{A} =
  \begin{bmatrix}
    \ddots & \vdots & \iddots \\
    \cdots & A_{s b} & \cdots \\
    \iddots & \vdots & \ddots
  \end{bmatrix}
  \in \RR^{S \times B}, \quad 
  \Vec{\psi}(\Vec{\phi}) =
  \begin{bmatrix}
    \vdots \\
    \psi_b(\Vec{\phi}) \\
    \vdots
  \end{bmatrix}
  \in \RR^B.
\end{equation}
The $\psi_b$'s are often chosen to be mononomials with maximum degree $M \geq 1$ in the state vector $\Vec{\phi}$; for example, with $M = 2$, we have
\begin{equation}
  \Vec{\psi}(\Vec{\phi}) =
  \begin{bmatrix}
    1 & \phi_1 & \cdots & \phi_S & \phi_1^2 & \cdots & \phi_S^2 & \phi_1 \phi_2 & \cdots & \phi_{S - 1} \phi_S
  \end{bmatrix}^\top.
\end{equation}
For general $M \geq 0$, we obtain a total $B = \binom{S + M}{S}$ monomials of degree up to $M$.
This number grows rapidly with increasing maximum degree for a fixed system size.
To avoid too many candidates contributing to the dynamics at once, SINDy constrains the unknown coefficients $\Mat{A}$ to be a sparse matrix.

The coefficients are learned in a data-driven manner: we collect the states and corresponding dynamics evaluations at discrete times $t_k$ along a system trajectory:
\begin{equation}
  \{(\Vec{\phi}_k, \Vec{r}_k) : \Vec{\phi}_k = \Vec{\phi}(t_k), \ 
                                \Vec{r}_k = \Vec{r}(\Vec{\phi}_k), \ 
        0 \leq k \leq K\},
\end{equation}
and use \eref{eq:sindy_factorization} to set up the linear system:
\begin{equation}
  \label{eq:sindy_linear_system}
  \underbrace{
    \begin{bmatrix}
        \Vec{r}_0 & \cdots & \Vec{r}_K
    \end{bmatrix}
  }_{\Mat{R} \in \RR^{S \times K}}
  \approx
  \Mat{A}
  \underbrace{
    \begin{bmatrix}
        \Vec{\psi}(\Vec{\phi}_0) & \cdots & \Vec{\psi}(\Vec{\phi}_K)
    \end{bmatrix}
  }_{\Mat{\Psi} \in \RR^{B \times K}}.
\end{equation}
We then use sparse regression to recover the coefficients; in particular, we employ the sequential thresholded least squares algorithm outlined in \cref{brunton2016discovering}:
starting with the ordinary least squares solution of \eref{eq:sindy_linear_system}, we disregard any basis $\psi_b$ with coefficients $A_{sb}$ below a specified threshold $\tau > 0$, and repeat the process with just the active set of basis functions.
We repeat this iteration until no new basis functions are eliminated, or a specified maximum number of elimination rounds is reached.
While setting up the linear system, a fixed percentage of the trajectory data is put aside for a validation study over the sparsity threshold.
As we increase the tolerance $\tau$, the sparsity of the coefficients $\Mat{A}$ increases along with the residual error $\| \Mat{R} - \Mat{A} \Mat{\Psi} \|$; along the Pareto front of this two-objective optimization (increasing sparsity vs.\ reducing error), we pick the tolerance that optimally balances the two. 
This data-driven coefficients discovery yields a surrogate model for the dynamics: $\Vec{\hat{r}}(\Vec{\phi}) = \Mat{A}\, \Vec{\psi}(\Vec{\phi})$.
We integrate this dynamics to obtain approximate system trajectories.

For forward uncertainty propagation in differential equations, the dynamics also depends on the model parameters $\Vec{\lambda}$. 
In the next section, we adapt the SINDy framework to construct surrogate dynamics for each sampled parameter vector by combining it with stochastic collocation.


%% file: figures/sparse_grid/demo.tex

\begin{figure}
  \centering
  \begin{subfigure}{0.23\textwidth}
    \centering
    \includegraphics[width=\linewidth]{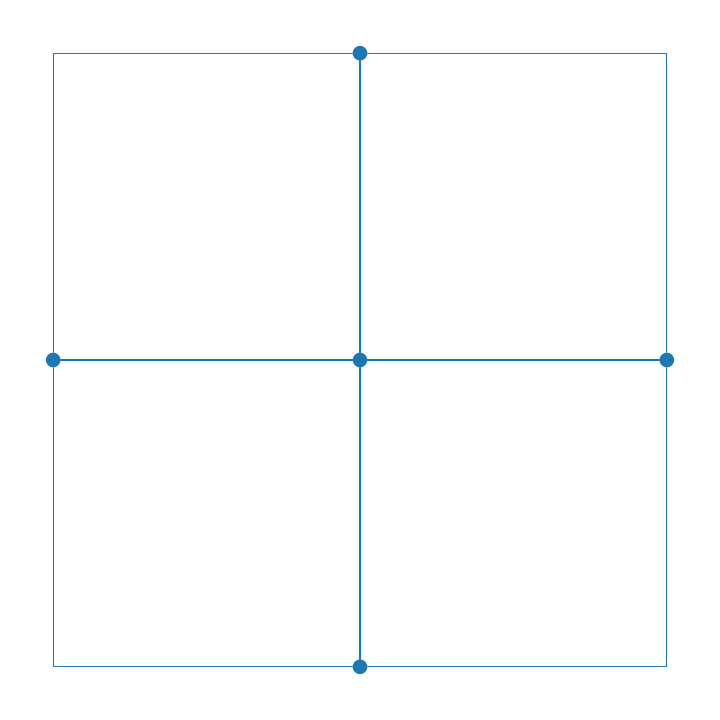}
    \caption{1st Order}
  \end{subfigure}
  ~
  \begin{subfigure}{0.23\textwidth}
    \centering
    \includegraphics[width=\linewidth]{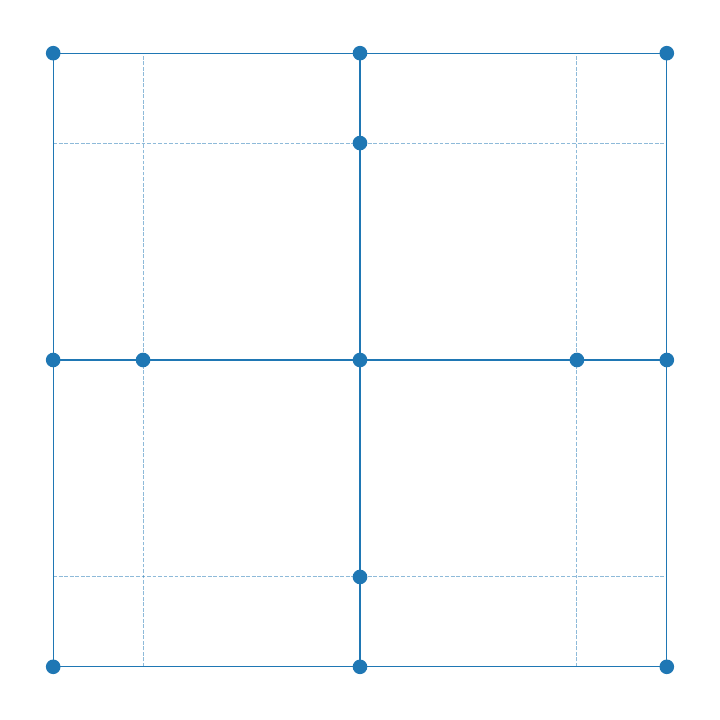}
    \caption{2nd Order}
  \end{subfigure}
  ~
  \begin{subfigure}{0.23\textwidth}
    \centering
    \includegraphics[width=\linewidth]{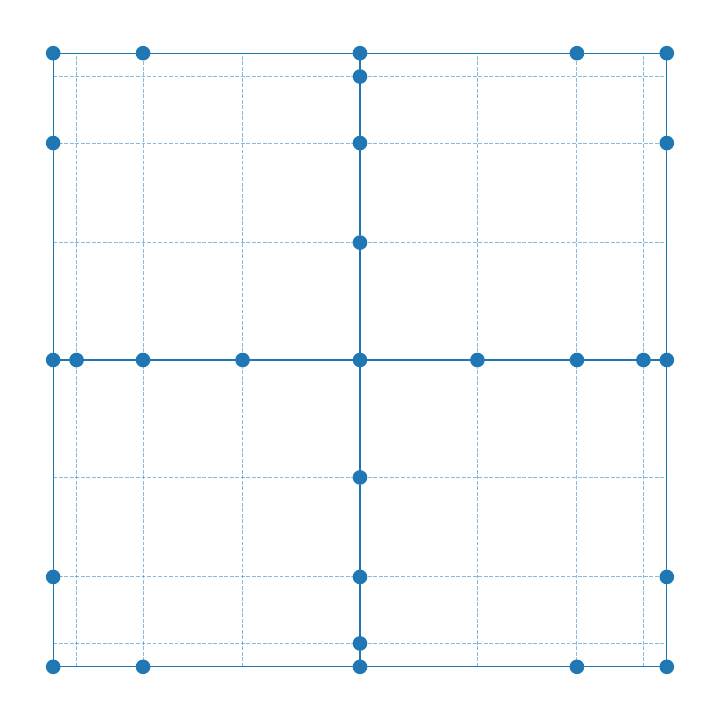}
    \caption{3rd Order}
  \end{subfigure}
  ~
  \begin{subfigure}{0.23\textwidth}
    \centering
    \includegraphics[width=\linewidth]{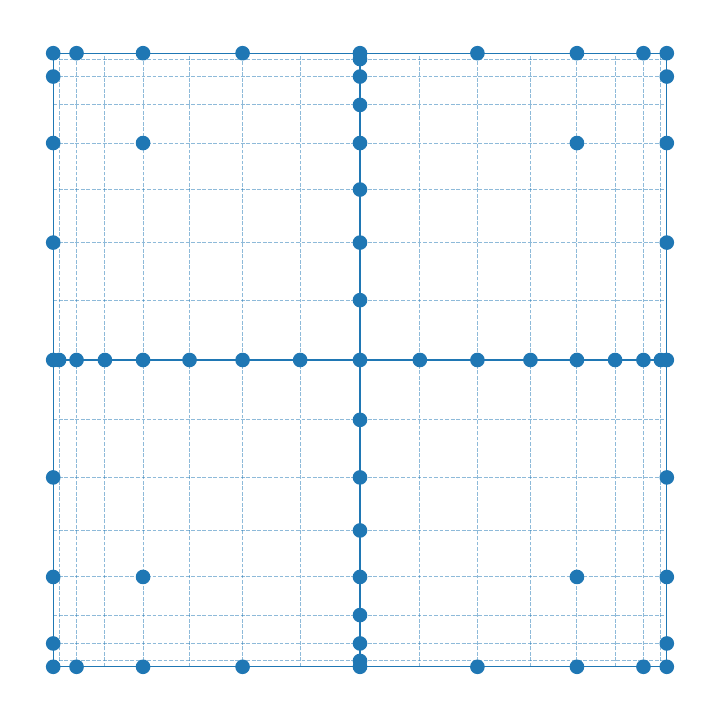}
    \caption{4th Order}
  \end{subfigure}
  \caption{
    Location of collocation nodes, marked by solid circles, using Clenshaw-Curtis sparse grids of various orders in the two-dimensional domain $[-1, 1]^2$.
    By systematically eliminating tensor grid nodes (corresponding to the `unmarked' crossings), sparse grids reduce the number of exact QoI evaluations.
    From left to right, there are 5, 13, 29, and 65 nodes in the sparse grids, and 9, 25, 81, and 289 nodes in the corresponding full grids.
  }
  \label{fig:sparse_grid_demo}
\end{figure}


%% file: sections/dynamics_sc.tex

\section{Stochastic Collocation over Dynamics} \label{sec:dynamics_sc}

Recall that, we aim to construct the state distribution $p_{\Vec{\phi}}(\Vec{\phi}; t)$ from system trajectories $\StateSmpExact{i}$ at Monte Carlo parameter samples $\ParamSmp{i}$.
To reduce the computational cost, we simulate exact system trajectories $\StateNdsExact{j}$ at collocation nodes $\ParamNds{j}$, and use \eref{eq:sc_state} to construct the surrogate trajectories $\StateSmpStdSC{i}$.
The accuracy of these State SC surrogates strongly depends on the complexity of the parameters-to-state map relative to the level of the Clenshaw-Curtis interpolation grid used to construct the collocation nodes.
With increasing non-linearity of the $\Vec{\lambda} \mapsto \Vec{\phi}(t, \Vec{\lambda})$ dependence, we need to use interpolation grids of increasingly higher levels, which correspondingly increases the number of exact simulations at the parameter collocation points.

Echoing \cref{jones2021minimally}, we maintain that directly approximating the dynamics with stochastic collocation (Dynamics SC) will reduce the number of exact simulations necessary for an accurate surrogate.
The parameter-to-dynamics map $\Vec{\lambda} \mapsto \Vec{r}(\Vec{\phi}, \Vec{\lambda})$ is typically less non-linear compared to the parameter-to-state map in ODE systems.
As an illustrative example, consider the one-dimensional, one-parameter ODE $\dot{\phi} = \lambda \phi$.
The dynamics $r = \lambda \phi$ depends linearly on the parameter, while the trajectory with initial condition $\phi(t = 0, \lambda) = 1$ is given by:
\begin{equation}
  \phi(t, \lambda) = \exp(\lambda t) = 1 + \frac{\lambda t}{1!} +
  \frac{\lambda^2 t^2}{2!} + \cdots,
\end{equation}
and it is highly non-linear in the parameter $\lambda$.
Approximating this trajectory as a polynomial will always introduce some error; however we can exactly capture the dynamics with two collocation points in the parameter space.
In this section, we leverage this observation, and combine SINDy with SC to construct surrogates that approximate the dynamics over both the state and parameters in a tractable manner.

\subsection{Dynamics Factorization}

Let us denote the dynamics corresponding to the collocation nodes and sample points in the parameter space, where we run exact simulations and construct approximations respectively, as $\DynamicsNdsExact{j} \equiv \Vec{r}(\Vec{\phi}, \ParamNds{j})$ and $\DynamicsSmpExact{i} \equiv \Vec{r}(\Vec{\phi}, \ParamSmp{i})$.
Stochastic collocation relates the two with an approximation of the form:
\begin{equation}
  \label{eq:sc_dynamics}
  \DynamicsSmpExact{i} \approx \DynamicsSmpDynSC{i} = \sum_{j = 1}^J \gamma_{ij} \, \DynamicsNdsExact{j}, \quad 1 \leq i \leq I.
\end{equation}
The `hat' in $\DynamicsSmpDynSC{i}$ again indicates an approximation, and the superscript `d-sc' indicates that stochastic collocation is carried out over the dynamics.
We then integrate this surrogate dynamics to construct approximate trajectories $\StateSmpDynSC{i}$.
We may use an explicit, quasi-static, or implicit time stepping scheme depending on the problem under consideration; but for the sake of simplicity, let us frame our method development in the forward Euler scheme:
\begin{equation}
  \StateSmpDynSC[t_{k + 1}]{i} = \StateSmpDynSC[t_k]{i} + \Delta t \, \DynamicsSmpDynSC[\StateSmpDynSC[t_k]{i}]{i}
\end{equation}
Note that, to use \eref{eq:sc_dynamics} directly in the RHS of this time evolution, we need to be able to evaluate the exact dynamics $\Vec{r}^\text{nds}_j$ corresponding to the \emph{collocation nodes} along the approximate trajectories $\StateSmpDynSC{i}$ which correspond to the \emph{parameter samples}.
These values are generally not computed during the exact simulations, which prevents us from using the dynamics approximation as is.

Driven by the need to evaluate surrogate dynamics at arbitrary state vector $\Vec{\phi}$ required for incorporation with ODE integrators, we propose a \emph{functional factorization of the dynamics} that separates the state variables from stochastic parameters:
\begin{equation}
  \label{eq:dynamics_factorization}
  \Vec{r}(\Vec{\phi}, \Vec{\lambda}) \approx \Mat{A}(\Vec{\lambda}) \Vec{\psi}(\Vec{\phi}).
\end{equation}
This form of the dynamics factorization is motivated by the SINDy framework; the novelty is that the sparse coefficients $\Mat{A}$ now depend on the value of the parameter.
The basis functions $\Vec{\psi}$, on the other hand, do not depend on the parameters; this point is crucial in rewriting \eref{eq:sc_dynamics} as:
\begin{equation}
  \DynamicsSmpDynSC{i}
  = \sum_{j = 1}^J \gamma_{ij} \, \DynCoeffNdsExact{j} \Vec{\psi}(\Vec{\phi}),
  \quad 1 \leq i \leq I.
\end{equation}
Here $\DynCoeffNdsExact{j}$ are the learned approximations of the coefficients $\Mat{A}(\ParamNds{j})$ satisfying $\DynamicsNdsExact{j} \approx \DynCoeffNdsExact{j} \Vec{\psi}(\Vec{\phi})$.
Note that we can rewrite this approximate dynamics as $\DynamicsSmpDynSC{i} = \DynCoeffSmpDynSC{i} \Vec{\psi}(\Vec{\phi})$ where
\begin{equation}
  \label{eq:sc_coeffs}
  \DynCoeffSmpDynSC{i} = \sum_{j = 1}^J \gamma_{ij} \, \DynCoeffNdsExact{j}, \quad 1 \leq i \leq I.
\end{equation}
This suggests the following workflow:
\begin{enumerate}
  \item
    Select a fixed set of basis functions $\Vec{\psi}$ for the system under consideration.
  \item
    Along each of the exactly simulated trajectories $\StateNdsExact{j}$ associated with the collocation nodes $\ParamNds{j}$, collect the model states \emph{and the dynamics evaluations} to set up the linear system \eref{eq:sindy_linear_system}. Then use sequentially thresholded least squares to learn the dynamics coefficients $\DynCoeffNdsExact{j}$ individually at the collocation nodes.
  \item
    Use \eref{eq:sc_coeffs} to compute the approximate dynamics coefficients $\DynCoeffSmpDynSC{i}$ associated with the sample points $\ParamSmp{i}$, and construct the dynamics surrogates $\DynamicsSmpDynSC{i}$.
  \item
    For each parameter sample, use surrogate dynamics with an appropriate ODE integrator to compute the approximate trajectories $\StateSmpDynSC{i}$.
\end{enumerate}
These approximate trajectories are significantly more accurate than the corresponding State SC trajectories $\StateSmpDynSC{i}$ as defined in \eref{eq:sc_state};  we demonstrate this with various numerical experiments in \sref{sec:results}.

\subsection{Time Adaptive Semi-Online Dynamics Learning}

Learning dynamics surrogates using SINDy requires trajectory data from previous time steps.
Thus, the Dynamics SC surrogates \emph{cannot} be constructed in a purely online fashion (at each time step without accessing prior history), unlike the State SC approach.
If storing the time stepping history is not a concern, the offline setup outlined in \pref{alg:dynsc_offline} can be used.

\input{algorithms/dynsc_offline}
\input{algorithms/dynsc_semi_online}

Alternatively, we note that while the response of a dynamical system can change drastically with variation in the initial conditions or system parameters, once these are fixed, the form of dynamics itself does not typically change over the course of the simulation.
In this setup, once we learn the dynamics from a trajectory segment of sufficient length, we can continue to use it for the remainder of the simulation, without needing to access the trajectory history.
This leads to a semi-online version of Dynamics SC, where:
\begin{itemize}
\item
  We simulate exact trajectories at the collocation nodes for a sufficient number of time steps.
  This number will be bounded above by how many instances of the full system state we can afford to save, across all the collocation nodes, within some allocated memory budget.
\item
  Then we compute the divergence between the exactly simulated trajectories at the collocation nodes, and the corresponding approximate trajectories predicted by the learned dynamics.
\item
  If the divergence is beyond some user-prescribed tolerances, we continue collecting data from the exact simulations and repeat the SINDy learning and validation.
  We stop this cycle when the learned dynamics have converged to the true dynamics, or we have exhausted the memory budget.
\item
  Once learning is complete, we compute the approximate dynamics at the sample parameters using SC, and continue generating approximate trajectories.
\item
  We consistently check the divergence criteria over the course of the simulation, and re-learn the dynamics as needed.
\end{itemize}
We present the salient details of this scheme in \pref{alg:dynsc_semi_online}.
Note that we use the forward Euler integrator in this pseudocode only to concretize and simplify the exposition; the method is amenable to other time stepping schemes (the results in \sref{sec:results} employ the fourth order Runge-Kutta algorithm for the Lorenz system, and the St\"{o}rmer-Verlet algorithm for solid mechanics PDE systems).

In this semi-online version, only lines 17 through 22 of \pref{alg:dynsc_semi_online}  collect and access histories of the exactly simulated trajectories; the rest of the steps rely at most on the immediately preceding time step.
Additionally, as the method gradually learns the complex dynamics in greater detail, we expect the divergence between estimated and simulated trajectories to occur less frequently.
By limiting the number of time steps used to construct the SINDy linear system in \eref{eq:sindy_linear_system}, we also reduce the computational cost of solving for the sparse coefficients.

\subsection{Localizing Dynamics Surrogates for PDE Systems}

The Dynamics SC framework, with SINDy-based dynamics factorization, works very well for low-dimensional ODE systems.
However, with increasing system dimensionality $S$ and fixed maximum degree $M$ of the monomial candidate functions, the total number of candidate functions grows very rapidly: $B = \binom{S + M}{S} \sim S^M$, and solving the resulting linear systems like \eref{eq:sindy_linear_system} becomes prohibitively costly.
This is especially true for ODE systems obtained by discretizing PDE systems on a computational mesh.

In practice, a majority of large-scale FEM implementations use linear basis functions with local support \cite{flanagan1981uniform}.
Continuing with notation from the solid mechanics setup presented in \sref{sec:ode_from_pde}, this choice implies that displacements $\Vec{u}_{n'}$ from only the immediate neighbors of node $n$ affect the computation of the forcing $\Vec{f}_n$.
In addition, the commonly employed lumped mass approximation renders the mass matrix $\Mat{M}$ diagonal; consequently the dynamics evaluations follow a similar pattern.
Since the number of neighbors of node $n$, denoted $N^\text{ngh}_n$, is usually much smaller than the number of nodes in the FEM mesh ($N^\text{ngh}_n \ll N$), this drastically reduces the number of the active variables in the SINDy approximation of the local dynamics $\Vec{r}_n$ from $S = D N$ to $S_n = D N^\text{ngh}_n$, with $S_n \ll S$.

Ideally, a properly trained SINDy model should automatically detect this sparsity structure from trajectory data; however:
\begin{itemize}
  \item
    We can easily determine the sparsity structure from the FEM mesh before any PDE simulations are initiated, and enforcing it directly on the SINDy model eliminates a potentially major source of numerical errors up-front.
  \item
    Reducing the number of active variables keeps the size of the linear systems in \eref{eq:sindy_linear_system} small, and we can solve these smaller system more efficiently.
In addition, the linear systems at each node, once constructed, are independent and can be solved in an embarrassingly parallel fashion.
\end{itemize}
Additionally, this approach focuses on discovering local behavior of the dynamics of the discretized PDE system independently at each node.
This is especially advantageous when the global dynamics is complex, e.g.\ in solid objects composed of multiple materials with different elastic parameters (see the notch example in \sref{sec:notch}).
By targeting a local approximation, we can use simple SINDy basis functions inferred from our knowledge of the constitutive equations.

We conclude our methodological developments by remarking that this local dynamics discovery strategy is not particular to the finite element method; we can apply it just as easily to finite difference or finite volume discretizations with local influence as well.
The corresponding improvement in efficiency makes the Dynamics SC framework particularly feasible for PDE systems.


%% file: algorithms/dynsc_offline.tex

\begin{algorithm}[t]
  \caption{Offline Dynamics SC}
  \label{alg:dynsc_offline}
  \begin{algorithmic}[1]
    \Require {Sample points $\ParamSmp{i}$, collocation nodes $\ParamNds{j}$, interpolation coefficients $\gamma_{ij}$, basis functions $\Vec{\psi}$}
    \Statex
    \For {$j = 1, \ldots, J$} \Comment{Loop over collocation nodes}
      \State {Collect state and dynamics evaluations $\StateNdsExact[t_k]{j}$ and $\DynamicsNdsExact[\StateNdsExact[t_k]{j}]{j}$ from exact simulation}
      \State {Form the linear system in \eref{eq:sindy_linear_system} and learn the coefficients $\DynCoeffNdsExact{j}$ using SINDy}
    \EndFor
    \For {$i = 1, \ldots, I$} \Comment{Loop over sample points}
      \State {Approximate the coefficients $\DynCoeffSmpDynSC{i}$ using the SC interpolation in \eref{eq:sc_coeffs}}
      \State {Integrate $\DynamicsSmpDynSC{i} = \DynCoeffSmpDynSC{i} \Vec{\psi}(\Vec{\phi})$ to obtain approximate trajectories $\StateSmpDynSC[t_k]{i}$}
    \EndFor
  \end{algorithmic}
\end{algorithm}


%% file: algorithms/dynsc_semi_online.tex

\begin{algorithm}[t]
  \caption{Semi-Online Dynamics SC (with Forward Euler Integrator)}
  \label{alg:dynsc_semi_online}
  \begin{algorithmic}[1]
    \Require {SC parameters $\ParamSmp{i}$, $\ParamNds{j}$, $\gamma_{ij}$; SINDy basis $\Vec{\psi}$; tolerances $\tau_\text{abs}$, $\tau_\text{rel}$}
    \Statex
    \State {Initialize time step index $k = 0$, coefficients $\DynCoeffNdsExact{j} = \Mat{0}$ for $1 \leq j \leq J$ and $\DynCoeffSmpDynSC{i} = \Mat{0}$ for $1 \leq i \leq I$}
    \While {$k \leq k_\text{max}$}
      \For {$j = 1, \ldots, J$} \Comment{Exact simulation over collocation nodes}
        \State {Advance using exact dynamics: $\StateNdsExact[t_{k + 1}]{j} = \StateNdsExact[t_k]{j} + \Delta t \; \DynamicsNdsExact[\StateNdsExact[t_k]{j}]{j}$}
        \State {Advance using learned dynamics: $\StateNdsApprx[t_{k + 1}]{j} = \StateNdsExact[t_k]{j} + \Delta t \; \DynCoeffNdsExact{j}  \; \Vec{\psi}(\StateNdsExact[t_k]{j})$}
        \If {$\| \StateNdsApprx[t_{k + 1}]{j} - \StateNdsExact[t_{k + 1}]{j} \| \geq \max\{\tau_\text{abs}, \tau_\text{rel} \| \StateNdsExact[t_{k + 1}]{j} \|\}$}
          \State {Flag learned dynamics as outdated}
        \EndIf
      \EndFor
      \If {Learned dynamics is not outdated} \Comment{Approximation over sample points}
        \For {$i = 1, \ldots, I$}
          \State {Advance surrogate trajectories: $\StateSmpDynSC[t_{k + 1}]{i} = \StateSmpDynSC[t_k]{i} + \Delta t \; \DynCoeffSmpDynSC{i} \; \Vec{\psi}(\StateSmpDynSC[t_k]{i})$}
        \EndFor
        \State {$k \gets k + 1$}
        \State {\textbf{go to} line 2}
      \EndIf
      \For {$j = 1, \ldots, J$} \Comment{Re-learn dynamics at collocation nodes}
        \For {$k' = 1, \ldots, K'$} \Comment{Collect trajectory segment of sufficient length $K'$}
          \State {Advance using exact dynamics: $\StateNdsExact[t_{k + k' + 1}]{j} = \StateNdsExact[t_{k + k'}]{j} + \Delta t \; \DynamicsNdsExact[\StateNdsExact[t_{k + k'}]{j}]{j}$}
        \EndFor
        \State {Form the linear system in \eref{eq:sindy_linear_system} and learn the coefficients $\DynCoeffNdsExact{j}$ using SINDy}
      \EndFor
      \For {$i = 1, \ldots, I$} \Comment{Re-interpolate dynamics at sample points}
        \State {Approximate the coefficients $\DynCoeffSmpDynSC{i}$ using the SC interpolation in \eref{eq:sc_coeffs}}
        \For {$k' = 1, \ldots, K'$} \Comment{Sync up with exact simulation time index}
          \State {Advance surrogate trajectories: $\StateSmpDynSC[t_{k + k' + 1}]{i} = \StateSmpDynSC[t_{k + k'}]{i} + \Delta t \; \DynCoeffSmpDynSC{i} \; \Vec{\psi}(\StateSmpDynSC[t_{k + k'}]{i})$}
        \EndFor
      \EndFor
      \State {$k \gets k + K' + 1$}
    \EndWhile
  \end{algorithmic}
\end{algorithm}


%% file: sections/results.tex

\section{Numerical Results} \label{sec:results}

We consider three model problems to test our methodology: the first is a system of ODEs that exhibit chaotic behavior, and the next two demonstrations are PDE systems from solid mechanics.
In each case, we treat system parameters as uncertain, while assuming that the initial conditions are known precisely.

\subsection{Error Metrics}

To recap our approach, we use the standard Monte Carlo ensemble to propagate uncertainty: we draw samples $\ParamSmp{i}$ from the parameter distribution $p_{\Vec{\lambda}}$ and run forward simulations to compute the sample trajectories $\StateSmpExact{i}$.
Using these samples, we construct the empirical distribution
\begin{equation}
  \label{eq:distribution_reference}
  \mathcal{D}^\text{smp}(t) = \left\{(\StateSmpExact{i}, 1 / I) : 1 \leq i \leq I\right\}
\end{equation}
for the target state distribution $p_{\Vec{\phi}}(\Vec{\phi}, t)$; here $1 / I$ is the equal weight assigned to each of the samples in the Monte Carlo ensemble.
We employ kernel density estimation (KDE) to visualize this distribution when required.
To reduce the cost of forward simulations of this brute-force approach, we employ State SC approximation and couple it with sparse grids to construct surrogate trajectories $\StateSmpStdSC{i}$.
We denote the corresponding state distribution as
\begin{equation}
  \mathcal{D}^\text{s-sc}(t) = \left\{(\StateSmpStdSC{i}, 1 / I) : 1 \leq i \leq I\right\}.
\end{equation}
This is the baseline against which we compare our proposed Dynamics SC approach, where we have the approximated trajectories $\StateSmpDynSC{i}$ and the corresponding state distribution
\begin{equation}
  \mathcal{D}^\text{d-sc}(t) = \left\{(\StateSmpDynSC{i}, 1 / I) : 1 \leq i \leq I\right\}.
\end{equation}

We use two metrics to compare the State and Dynamics SC approximations against the reference constructed using exact simulations. First, we directly compare the system trajectories using average and maximum errors over the samples:
\begin{align}
  \label{eq:error_metric_traj_avg}
  \operatorname{error}_\text{avg}^\text{method}(t) &= \frac{\frac{1}{I} \sum_{i = 1}^I \| \StateSmpExact{i} - \hat{\Vec{\phi}}^\text{method}_i(t) \|}{\phi_\text{max}}, \; \text{and} \\
  \label{eq:error_metric_traj_max}
  \operatorname{error}_\text{max}^\text{method}(t) &= \frac{\max_{1 \leq i \leq I} \| \StateSmpExact{i} - \hat{\Vec{\phi}}^\text{method}_i(t) \|}{\phi_\text{max}},
\end{align}
for $\text{method} \in \{\text{s-sc}, \text{d-sc}\}$. The scaling factor $\phi_\text{max}$ is the maximum magnitude of the state variable over all exactly simulated sample trajectories and all time:
\begin{equation}
  \phi_\text{max} = \max_{i, t} \| \StateSmpExact{i} \|;
\end{equation}
we use it to uniformly scale errors in the various state variables (e.g. displacement and velocity) of wildly different magnitudes, and construct a relative error estimate.

While these sample-to-sample metrics are standard indicators of surrogate performance, we are more interested in directly comparing the reference distribution $\mathcal{D}^\text{smp}(t)$, constructed from exact simulations at the Monte Carlo samples, against the distributions $\mathcal{D}^\text{method}(t)$ constructed with the SC surrogates.
For this, we use the 1-Wasserstein metric, which evaluates distance between two distributions defined over a metric space by solving an optimal transport problem \cite{villani2009optimal}.
In particular, for two normalized empirical distributions $\mathcal{D} = \{(\Vec{\phi}_i, w_i) : 1 \leq i \leq I\}$ and $\mathcal{D}' = \{(\Vec{\phi}'_j, w'_j) : 1 \leq j \leq J\}$, where $w_i$ and $w_j'$ represent the weights of the corresponding samples $\Vec{\phi}_i$ and $\Vec{\phi}'_j$, the Wasserstein distance is obtained by solving the linear program:
\begin{equation}
  \| \mathcal{D} - \mathcal{D}' \| := \min \bigg\{\sum_{i = 1}^I \sum_{j = 1}^J f_{i j} d(\Vec{\phi}_i, \Vec{\phi}'_j) \; \bigg\vert \; f_{i j} \geq 0, \; \sum_{j = 1}^J f_{i j} \leq w_i, \; \sum_{i = 1}^I f_{i j} \leq w'_j, \; \sum_{i = 1}^I \sum_{j = 1}^J f_{i j} = 1\bigg\},
\end{equation}
where $d(\cdot, \cdot)$ is the distance metric of the underlying space \cite{komiske2019metric}.
We use the open-source \texttt{wasserstein}\footnote{https://github.com/pkomiske/Wasserstein} Python package to facilitate this computation.

\subsection{Lorenz System} \label{sec:lorenz}

The Lorenz system is a well-known three-dimensional nonlinear dynamical system with three parameters \cite{strogatz2018nonlinear}:
\begin{equation}
  \dot{x} = \sigma (y - x), \quad \dot{y} = x (\rho - z) - y, \quad \dot{z} = xy - \beta z.
\end{equation}
The parameter values of $\rho = 28$, $\sigma = 10$, and $\rho = \frac{8}{3}$ correspond to the characteristic ``butterfly wings''-shaped chaotic Lorenz attractor.
We treat each of these parameters as uncertain, and assume that they have been determined with up to 5\% error w.r.t.\ the nominal values:
\begin{equation}
  \rho \sim \Uniform\left[28 \pm \frac{7}{5}\right], \quad \sigma \sim \Uniform\left[10 \pm \frac{1}{2}\right], \quad \beta \sim \Uniform\left[\frac{8}{3} \pm \frac{2}{15}\right].
\end{equation}
To construct a reference distribution over the state $\Vec{\phi} = (x, y, z)$, we sample 512 parameter triples from this distribution over $\Vec{\lambda} = (\rho, \sigma, \beta)$.
Then, we simulate the corresponding system trajectories, starting from $\Vec{\phi}_0 = (1, 1, 1)$, using the fourth order Runge-Kutta integrator with time step size $\Delta t = 10^{-6}$ for $10^8$ steps (total simulation time of $t_\text{final} = 100$).
At any point during the time evolution, these trajectories serve as the reference empirical distribution of the state.
In \fref{fig:lorenz_distribution_demo}, we use kernel density estimation (KDE) to construct continuous representations of one- and two-dimensional projections of this empirical distribution at different stages of the simulation.
As can be seen from the marginal distributions, the state distribution is complex, with multimodal behavior that evolves with time.
We use these as reference distributions to evaluate the performance of the stochastic collocation surrogates.

\input{figures/lorenz/distribution_demo}

For a direct comparison of trajectory accuracy using the State and Dynamics SC frameworks, we compute approximate trajectories at the previously constructed 512 parameter sample points.
In both cases, we use 1st, 2nd, and 3rd order sparse grids in the three-dimensional parameter space, corresponding to 7, 25, and 67 collocation nodes where we run exact system simulations, to construct the approximations.
For the SINDy learning during Dynamics SC, we specify using a maximum of 500 dynamics evaluations, collected every $10^4$ integrator steps at $\Delta t = 10^{-2}$ intervals, to learn the forcing term with a 5th order polynomial expansion.
We noticed that 200 such evaluations (corresponding to initial 2\% of the full trajectories at each collocation node) were sufficient to learn the dynamics, and reuse it for the remaining 98\% of the simulations.

\input{figures/lorenz/trajectory_demo}

\input{figures/lorenz/error}

We plot a sample trajectory of the Lorenz system, constructed using the three methods (exact simulation, State SC, and Dynamics SC), in \fref{fig:lorenz_trajectory_demo} for visual comparison.
Clearly Dynamics SC produces a better approximation.

For a more quantitative comparison, we compute the deviations of the approximate trajectories from the corresponding exact simulation.
\fref{fig:lorenz_trajectory_error} tracks the mean errors, as defined in \eref{eq:error_metric_traj_avg}, for the two family of surrogates over time.
In all cases, the error grows over the course of simulation until it reaches the maximum\footnote{The maximum exits due to the bounded nature of the Lorenz system.}, which is characteristic of a chaotic system.
However, we observe that the rate of growth for the Dynamics SC error is significantly lower compared to State SC.

To compare the quality of the constructed state distributions, we use the Euclidean metric in the 3D state space $\Vec{\phi} = (x, y, z)$ to compute the 1-Wasserstein distances $\| \mathcal{D}^\text{smp}(t) - \mathcal{D}^\text{s-sc}(t) \|$ and $\| \mathcal{D}^\text{smp}(t) - \mathcal{D}^\text{d-sc}(t) \|$ between the distributions constructed from exact and approximate system trajectories.
In \fref{fig:lorenz_distribution_error}, we plot the evolution of these two errors; we again see that the Dynamics SC significantly outperforms the State SC.
Note that, there is a persistent gap between the distribution errors obtained from the two SC approaches, while there is no such discrepancy in the trajectory errors in \fref{fig:lorenz_trajectory_error}.
This is because comparison in distributions is inherently weaker (more forgiving): the Wasserstein distance is an aggregate measure over all the trajectories, and not a one-to-one matching over individual trajectories.

\subsection{One-Dimensional Bar Impact} \label{sec:bar}

Next, we consider a one-dimensional problem, where a solid bar of length $L$ is initially moving at a constant velocity.
At time $t = 0$, the right end of the bar impacts a fixed wall and is immediately brought to rest.
In this setup, the momentum balance in \eref{eq:solid} simplifies to
\begin{equation}
  \rho_0 u_{tt}(t, X) = \sigma_X(t, X), \quad t > 0, \quad 0 < X < L.
\end{equation}
where $\sigma$ denotes the stress. The boundary conditions at the two ends are given by:
\begin{equation}
  u_X(t, X = 0) = 0, \quad u(t, X = L) = 0,
\end{equation}
and the initial conditions are specified by:
\begin{equation}
  u(t = 0, X) = 0, \quad u_t(t = 0, X) =
  \begin{cases}
    v_0 & \text{if} \quad 0 \leq X < L, \\
    0 & \text{if} \quad X = L.
  \end{cases}
\end{equation}
We assume a nonlinear constitutive model for the stress:
\begin{equation}
  \label{eq:bar_impact_stress}
  \sigma(t, X) = \epsilon_1 u_X(t, X) (1 + \epsilon_2 u_X(t, X)^2);
\end{equation}
this results in a non-linear wave equation for the displacement field.
In our simulations, we discretize the bar using 100 uniform linear elements and employ a lumped mass approximation.
The following parameters are fixed: linear density $\rho_0 = 1$, bar length $L = 2$, initial velocity $v_0 = 1$.
The stress model parameters are stochastic with distributions $\epsilon_1, \epsilon_2 \sim \Uniform[1 \pm 0.05]$.
We draw 128 samples from the two-dimensional parameter distribution to construct the reference trajectories using the explicit St\"{o}rmer-Verlet integration with step size $\delta t = 10^{-3}$ for 20000 steps (total simulation time $t_\text{final} = 20$).

\input{figures/bar_impact/error_u}

\input{figures/bar_impact/error_v}

\input{figures/bar_impact/error_emd}

To compare the performance of the two SC approaches, we construct approximate trajectories at the sampled parameter points.
In both cases, we use 1st, 2nd, and 3rd order sparse grids, corresponding to 5, 13, and 29 exact simulations at the collocation points, to construct the surrogates.
In the Dynamics SC approach, we specify using a maximum number of 200 dynamics evaluations, at $\Delta t = 0.02$ time intervals, to learn the dynamics; however 100 evaluations (corresponding to initial 10\% of the simulation) proved sufficient for this purpose.
Once again, we used 5th degree polynomials to approximate the dynamics with SINDy.
In \fref{fig:bar_impact_error_u}, we plot the average and maximum errors in the displacement field $u(t, X)$, computed over the 128 parameter samples at each time step $t_k$ of the ODE integration and each node $X_n$ of the spatial discretization.
We use a logarithm scale-based color map to capture the variation in the error across eight orders of magnitude.
The impact occurs at the upper-left corner ($t = 0$, $X = 2$) of each of the panels, and the resulting shock wave propagates to the rest of the bar with increasing time.
Due to the finite wave speed, we see a triangular zero-error region on the lower-left of the panels, where the shock is yet to arrive.
The wavefront is scattered by the boundaries; we see artifacts of this in the error plots in the oblique boundaries between the error levels.
Note the fixed left boundary creates a thin region of zero error for both methods.
Comparing the two SC approaches, we observe a drastic improvement in the accuracy of Dynamics SC surrogates over State SC, for the 2nd and 3rd order sparse grids.
The 1st order sparse grid is an outlier in this respect; we do see improvements in the error, but they are not as drastic.
This occurs because the 1st order sparse grid fails to accurately evaluate the cross term ($\epsilon_1 \epsilon_2$) in the stress expression from \eref{eq:bar_impact_stress}.

We note the change in error statistics of the State SC as we increase the quality of the sparse grid: high-error regions are `pushed away' from the upper-left corner of the panels (where the shock originates).
This indicates that initially State SC surrogates learn better approximations with higher-level sparse grids, but that improvement is gradually lost as the simulation progresses.
Contrast this with the behavior of Dynamics SC: the 3rd level sparse grid does not yield any significant improvements over the 2nd level grid as the dynamics is sufficiently resolved with the lower-level grid.

The second primary solution from our simulations, velocity $\dot{u}(t, X)$ field, exhibit similar error characteristics.
We show the point-wise (at every time step and at each spatial node) average and maximum errors obtained from the two SC methods using the 2nd order sparse grid.
The overall behavior is nearly identical to that of the displacement field apart from one caveat: the errors are higher by an order of magnitude.
This loss of accuracy can be attributed to the velocity field being less smooth compared to displacement.

Next, we compare the distribution errors.
For each time step $t_n$ during the ODE integration and at each spatial node $X_n$, we construct three empirical distributions $\mathcal{D}^\text{smp}(t_k, X_n)$, $\mathcal{D}^\text{s-sc}(t_k, X_n)$, and $\mathcal{D}^\text{d-sc}(t_k, X_n)$ just like the Lorenz system demonstration; we use $\Vec{\phi}(t, x) = (u(t, x), \dot{u}(t, x))$ as the system state.
Using the Euclidean metric to compute the Wasserstein distances between these local distributions, we plot the errors, obtained from the two SC methods using the 1st and 2nd order sparse grids, in \fref{fig:bar_impact_error_emd}.  
We observe that the error statistics follow the already established patterns: (a) until the initial shock propagates through the spatial domain, both SC methods produce zero errors, (b) Dynamics SC produces lower errors compared to State SC, (c) the improvement is especially prevalent with the 2nd order sparse grid, which fully resolves the dynamics.

\subsection{Notch Impact} \label{sec:notch}

In this more complex quasi-3D example, a glass layer is bonded to an aluminum layer, and the aluminum layer has three notches that induce stress concentrations (see \fref{fig:notched_plate_setup}). 
For both materials we adopt a compressible, finite elastic neo-Hookean model
\begin{equation}
  \label{eq:hyperelastic}
  \Mat{P} \Mat{F}^\top = \Mat{F} \Mat{P}^\top =
  \frac{1}{2} \kappa \, (J^2 -1) \Mat{I} + \mu J^{-2/3} \left( \Mat{B} - \frac{1}{3} \tr(\Mat{B}) \Mat{I} \right),
\end{equation}
where $\Mat{B} = \Mat{F} \Mat{F}^T$ is the left Cauchy-Green deformation tensor, $\Mat{F} = \nabla_{\Vec{X}} \Vec{x}(t, \Vec{X})$ is the deformation gradient, and $J = \det\Mat{F}$ is the Jacobian.
The reference density $\rho_0$ of both materials are fixed, while the elastic moduli (bulk modulus $\kappa$ and shear modulus $\mu$) are treated as uncertain:
\begin{align*}
  \text{glass:}    && \rho_0 &= 2.7~\text{g/cm}^3, & \kappa &\sim \Uniform[30~\text{GPa}, 48~\text{GPa}], & \mu &\sim \Uniform[20~\text{GPa}, 34~\text{GPa}], \\
  \text{aluminum:} && \rho_0 &= 3.0~\text{g/cm}^3, & \kappa &\sim \Uniform[68~\text{GPa}, 70~\text{GPa}], & \mu &\sim \Uniform[25~\text{GPa}, 27~\text{GPa}].
\end{align*}
We use 512 samples from this four-dimensional stochastic parameter space to construct the state distributions.

\input{figures/notched_plate/setup}

The composite plate has fixed boundary conditions at the left end, and it is subject to an off center impact loading (indicated in \fref{fig:notched_plate_setup} with the black wall and grey arrows).
This impact induces an initial velocity of 40~m/s at the nodes along the bottom right boundary segment.
We discretize this structure with 2558 nodes and 1169 elements, and use the open-source Lagrangian solid mechanics code \texttt{NimbleSM}\footnote{https://github.com/NimbleSM/NimbleSM} to simulate system trajectories with the explicit Störmer-Verlet integrator for $10^4$ steps with step size $\Delta t = 10^{-9}$~s (total duration of $t_\text{final} = 0.01$~ms).
We employ the two SC frameworks to construct approximate trajectories at the sample parameters and compare them to the results from the exact simulations.
We use a first order sparse grid, corresponding to 9 exact simulations at the collocation points in the four-dimensional parameter space, to construct the surrogates.
For Dynamics SC, we used 1001 dynamics evaluations at $\Delta t = 10^{-8}$~s intervals over the entire course of the simulations to learn the node-local dynamics using linear SINDy basis.

\input{figures/notched_plate/error_u}

As with the previous demonstration, we first plot the average and maximum error magnitudes in the estimated displacement field  $\Vec{u}(t, \Vec{X})$, computed over the 512 sample trajectories at different times $t$ during the course of the simulation, in \fref{fig:notched_plate_error_u}.
We again use a logarithmic color scale to capture the wide range of error values.
The pattern created by reflections at various material boundaries, including the 3 notches, are again evident in these plots.
We note, again, that the error is always zero at the fixed left end of the plate.
During the initial phase of the simulation, State SC produces lower error compared to Dynamics SC where the shock has not yet arrived, e.g.\ the left-most blue `bulb' at $t = 2.5~\mu$s; the situation is reversed in the right-most red bulb, where the shock originates.
Any accuracy advantage of State SC is quickly lost as the simulation progresses and the wave propagates through the entire computational domain.
In \fref{fig:notched_plate_error_v}, we create a similar plot the velocity field $\dot{\Vec{u}}(t, \Vec{X})$ at the final instance  of the simulation ($t = 10~\mu$s), and show analogous dramatic improvement in the error when using Dynamics SC.

\input{figures/notched_plate/error_v}
\input{figures/notched_plate/error_emd}

Finally, we construct the three empirical distributions $\mathcal{D}^\text{smp}(t_k, \Vec{X}_n)$, $\mathcal{D}^\text{s-sc}(t_k, \Vec{X}_n)$, and $\mathcal{D}^\text{d-sc}(t_k, \Vec{X}_n)$ at each time step $t_n$ and FEM node $\Vec{X}_n$ over the six-dimensional states
\begin{equation}
  \Vec{\phi}_{\Vec{X}_n}(t_k) =
  \begin{bmatrix}
    \Vec{u}(t_k, \Vec{X}_n) / s_u \\
    \dot{\Vec{u}}(t_k, \Vec{X}_n) / s_{\dot{u}}
  \end{bmatrix},
\end{equation}
The factors $s_u = 10^{-3}$~cm and $s_{\dot{u}} = 10^3$~cm/s non-dimensionalize the displacement and velocity fields, and also bring their values to comparable orders of magnitudes.
These scaling constants are determined empirically, based on the magnitude of the displacement and velocity fields from the exact simulations.
We use the Euclidean metric in this six-dimensional state space to compute the Wasserstein distances {$\| \mathcal{D}^\text{smp} - \mathcal{D}^\text{s-sc} \|$} and {$\| \mathcal{D}^\text{smp} - \mathcal{D}^\text{d-sc} \|$} between the approximate and exactly simulated state distributions.
\fref{fig:notched_plate_error_emd} shows the evolution of this distribution error over the course of the simulation.
We observe features similar to those present in the primary displacement and velocity fields:
(a) in regions where the shock front is yet to arrive, State SC produces smaller error than Dynamics SC,
(b) as the simulation progresses, Dynamics SC eventually outperforms State SC in terms of reconstruction errors by orders of magnitude.


%% file: figures/lorenz/distribution_demo.tex

\begin{figure*}
  \centering

  \begin{subfigure}{0.48\textwidth}
    \centering
    \includegraphics[width=\linewidth]{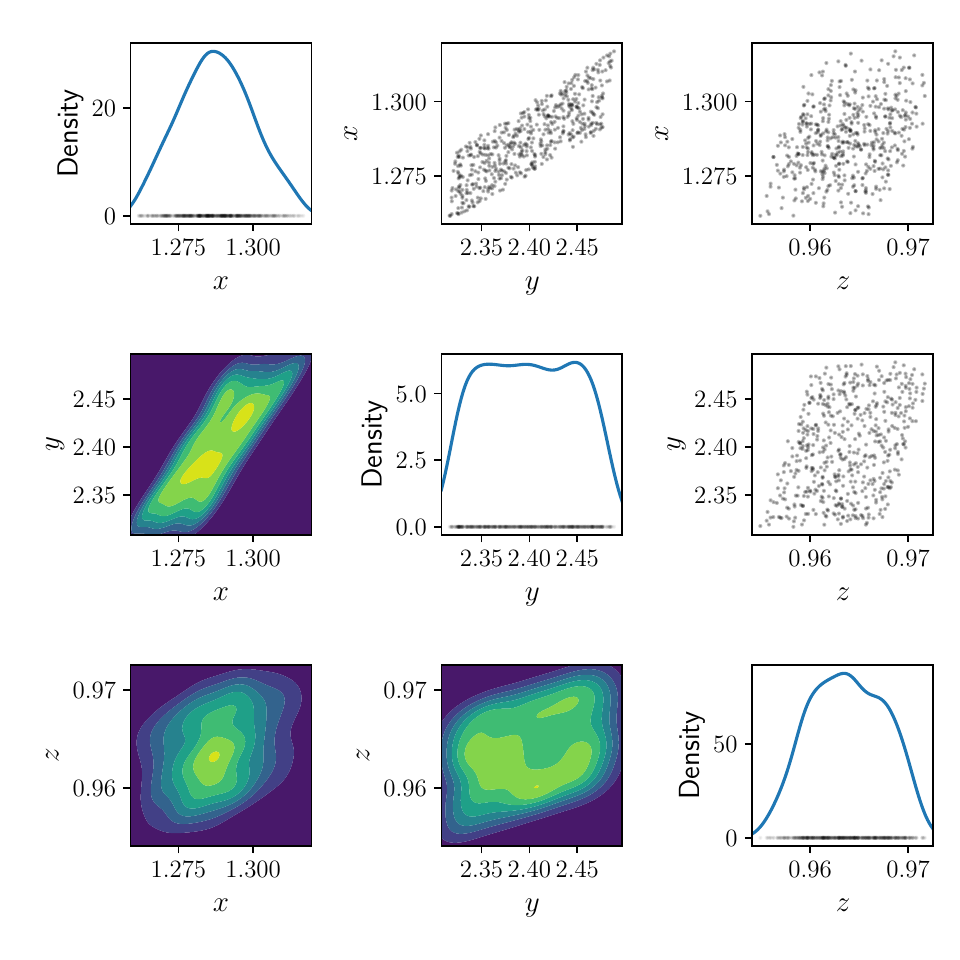}
    \caption{$t = 0.05$}
  \end{subfigure}
  ~
  \begin{subfigure}{0.48\textwidth}
    \centering
    \includegraphics[width=\linewidth]{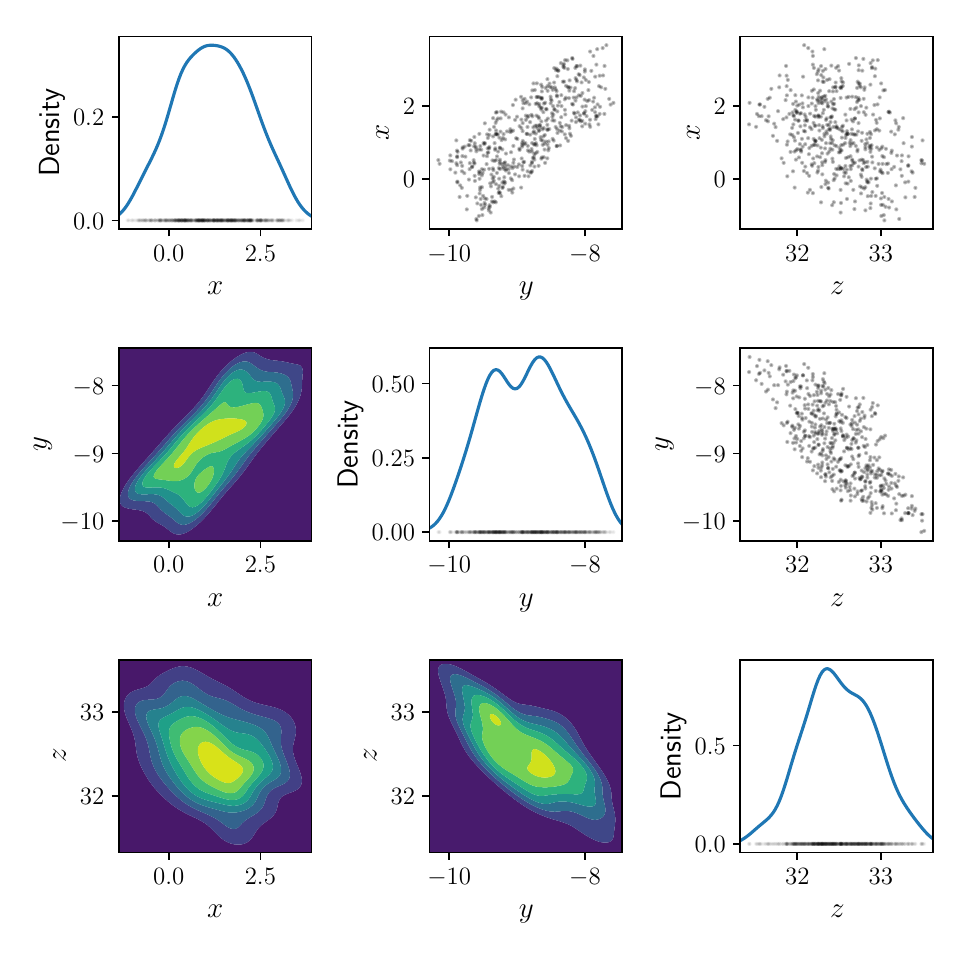}
    \caption{$t = 0.5$}
  \end{subfigure}

  \bigskip

  \begin{subfigure}{0.48\textwidth}
    \centering
    \includegraphics[width=\linewidth]{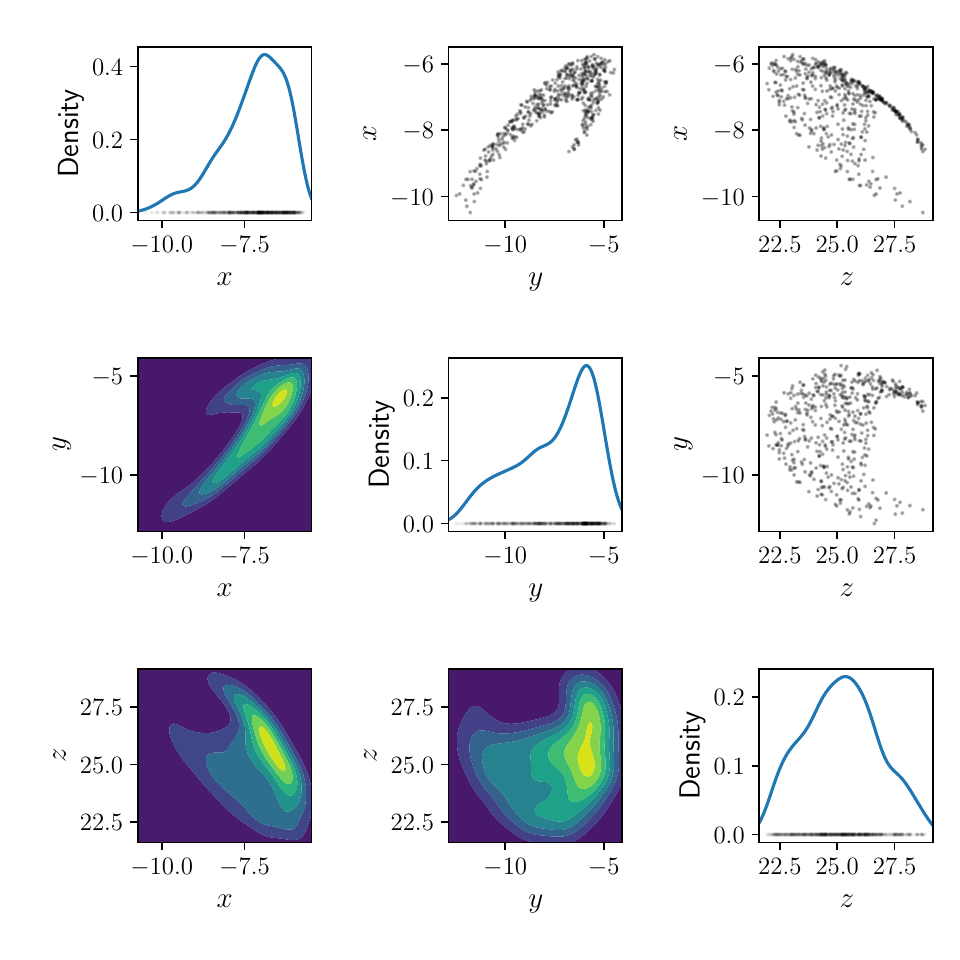}
    \caption{$t = 5$}
  \end{subfigure}
  ~
  \begin{subfigure}{0.48\textwidth}
    \centering
    \includegraphics[width=\linewidth]{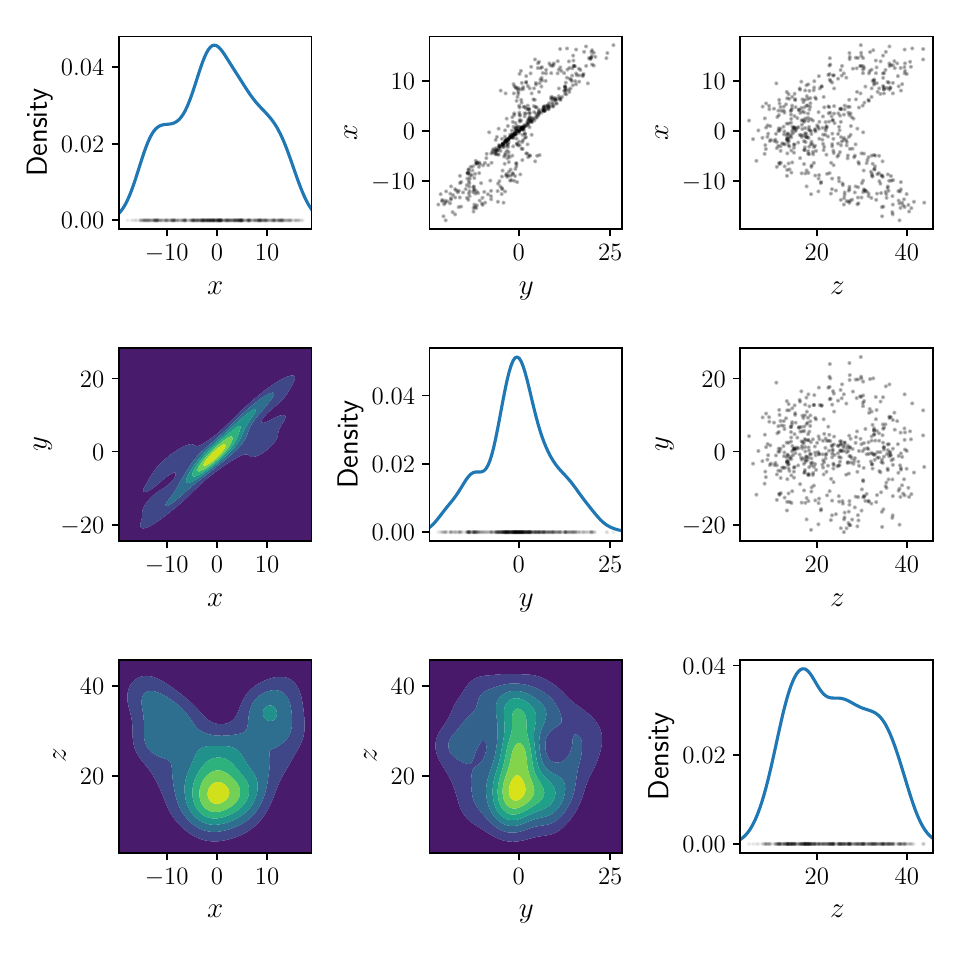}
    \caption{$t = 50$}
  \end{subfigure}

  \caption{
    State distributions $p_{x y z}(x, y, z; t)$ of the Lorenz system at various times, constructed by time stepping the exact ODE, projected onto one-dimensional coordinate axes (diagonal panels) and two-dimensional coordinate planes (lower triangular panels).
    The distributions are constructed using kernel density estimation (KDE) from 512 trajectories, each staring from $x_0 = y_0 = z_0 = 1$, and corresponding to individual $(\rho, \sigma, \beta)$ parameter samples.
    The state samples are shown in the scatter plots with black markers (diagonal and upper triangular panels).
  }
  \label{fig:lorenz_distribution_demo}
\end{figure*}


%% file: figures/lorenz/trajectory_demo.tex

\begin{figure}
  \centering
  \begin{subfigure}{0.32\textwidth}
    \centering
    \includegraphics[width=\linewidth]{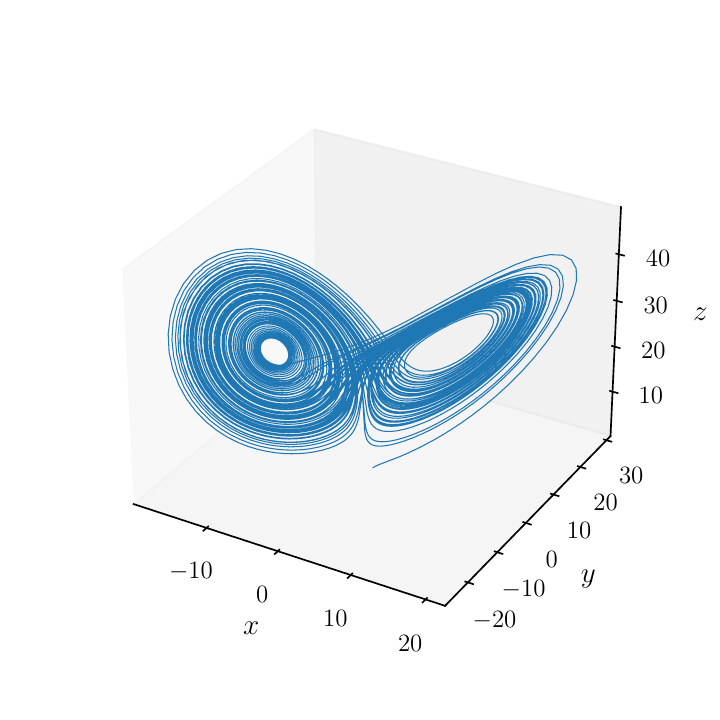}
    \caption{Exact}
  \end{subfigure}
  ~
  \begin{subfigure}{0.32\textwidth}
    \centering
    \includegraphics[width=\linewidth]{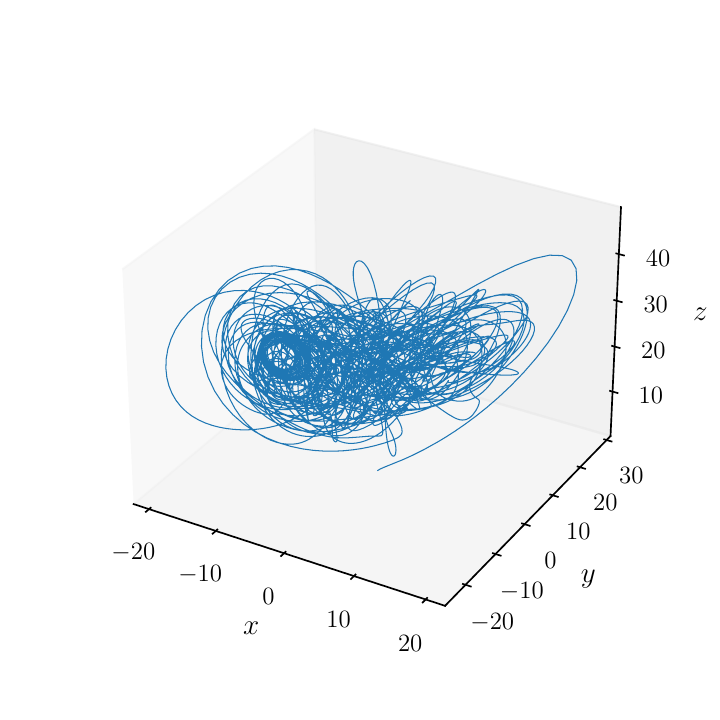}
    \caption{State SC}
  \end{subfigure}
  ~
  \begin{subfigure}{0.32\textwidth}
    \centering
    \includegraphics[width=\linewidth]{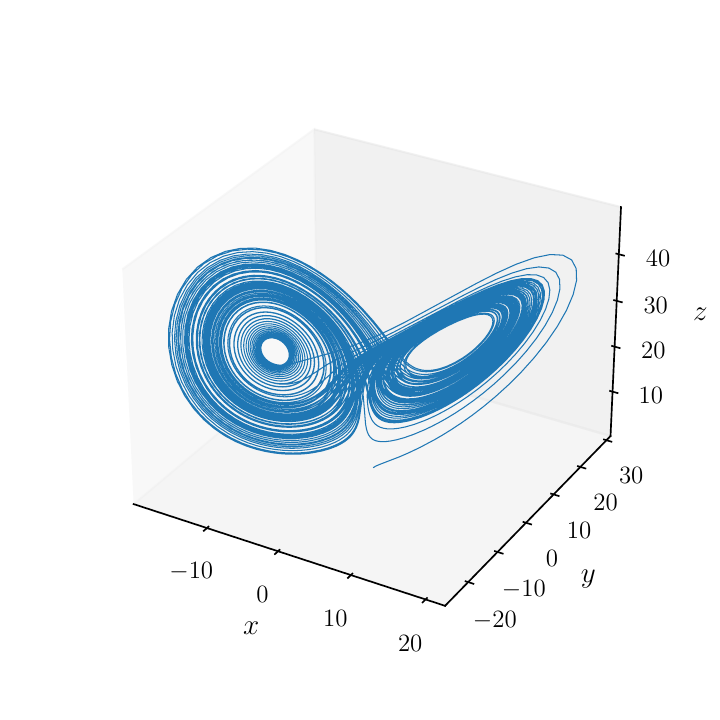}
    \caption{Dynamics SC}
  \end{subfigure}
  \caption{
    Visual comparison of approximate trajectories generated from the two SC approximations against the corresponding exact trajectory of the Lorenz system.
    We use the 1st order sparse grid (7 exact simulations) to construct these approximations.
    Our proposed Dynamics SC method clearly produces a better approximation of the system trajectory.
  }
  \label{fig:lorenz_trajectory_demo}
\end{figure}


%% file: figures/lorenz/error.tex

\begin{figure}
  \centering
  \begin{subfigure}{0.48\textwidth}
    \centering
    \includegraphics[width=\linewidth]{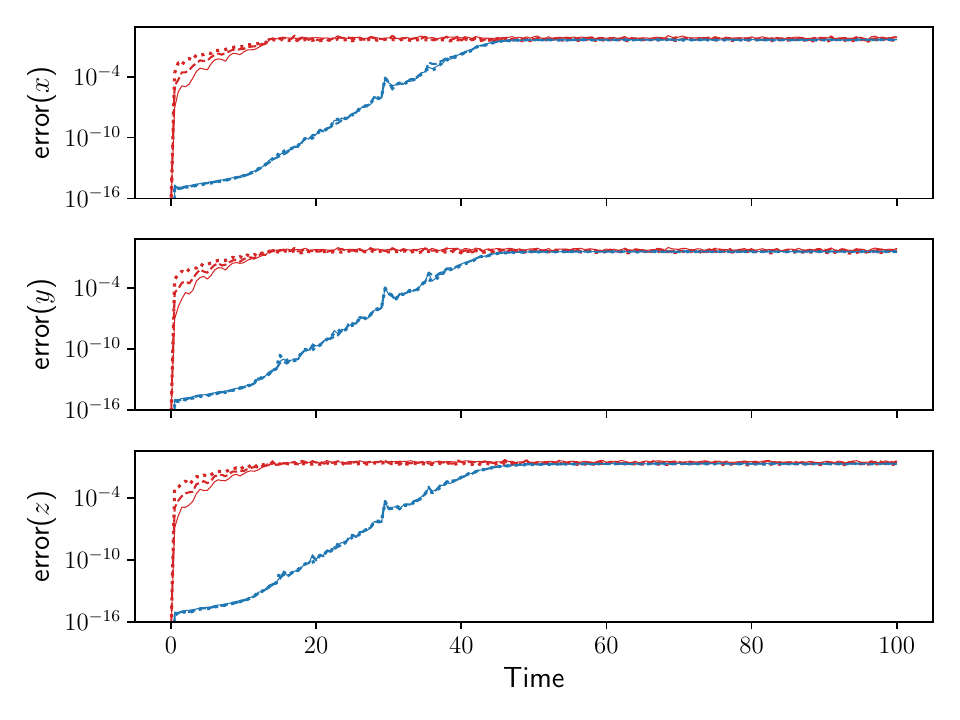}
    \caption{Trajectory Error}
    \label{fig:lorenz_trajectory_error}
  \end{subfigure}
  ~
  \begin{subfigure}{0.48\textwidth}
    \centering
    \includegraphics[width=\linewidth]{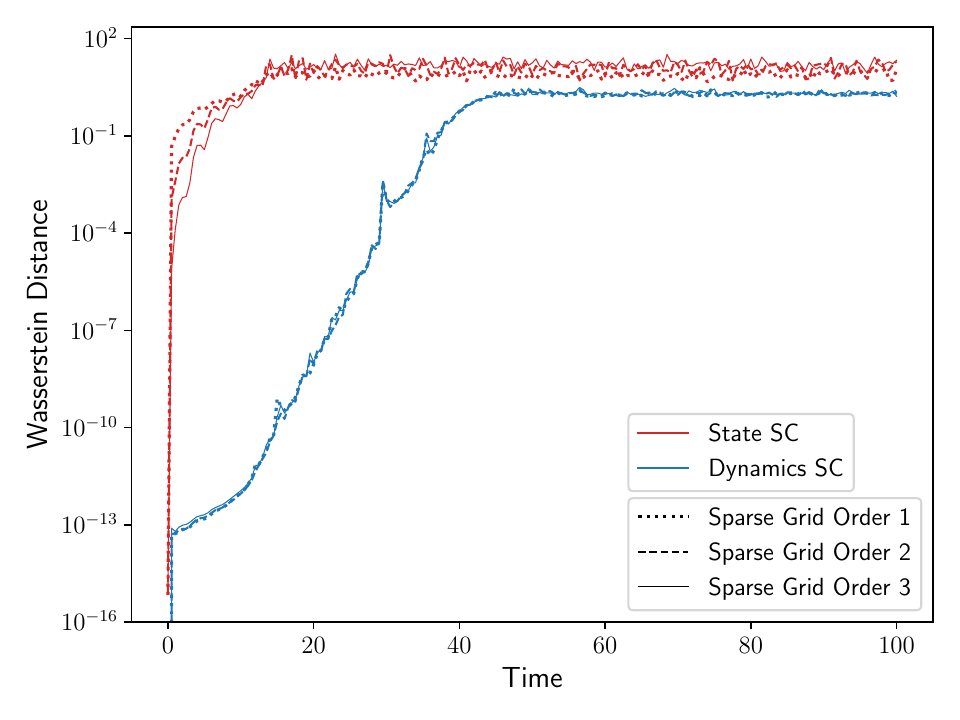}
    \caption{Distribution Error}
    \label{fig:lorenz_distribution_error}
  \end{subfigure}
  \caption{
    Time evolution of various metrics comparing the accuracy of stochastic collocation surrogates against the reference exact simulations in the Lorenz system.
    The left panel shows the deviations of the three components of the approximate state from the simulated exact values, averaged over the 512 parameter samples.
    The right panel captures the evolution of the Wasserstein distance of the empirical state distributions from the reference.
    In each case, we plot the errors corresponding to the State and Dynamics SC surrogates (in blue and red colors, respectively) constructed with first, second, and third order sparse grid collocation points (in dotted, dashed and solid curves, respectively).
    As expected, increasing the number of exact simulations (by using higher-order sparse grid) does improve the accuracy of the State SC surrogates.
    But, even with 67 exact simulations (3rd order grid), the accuracy of these surrogates is much lower than the Dynamics SC surrogates constructed from 7 exact simulations (1st order grid).
  }
  \label{fig:lorenz_error}
\end{figure}


%% file: figures/bar_impact/error_u.tex

\begin{figure*}
  \centering

  \begin{subfigure}{0.48\textwidth}
    \centering
    \includegraphics[width=\linewidth]{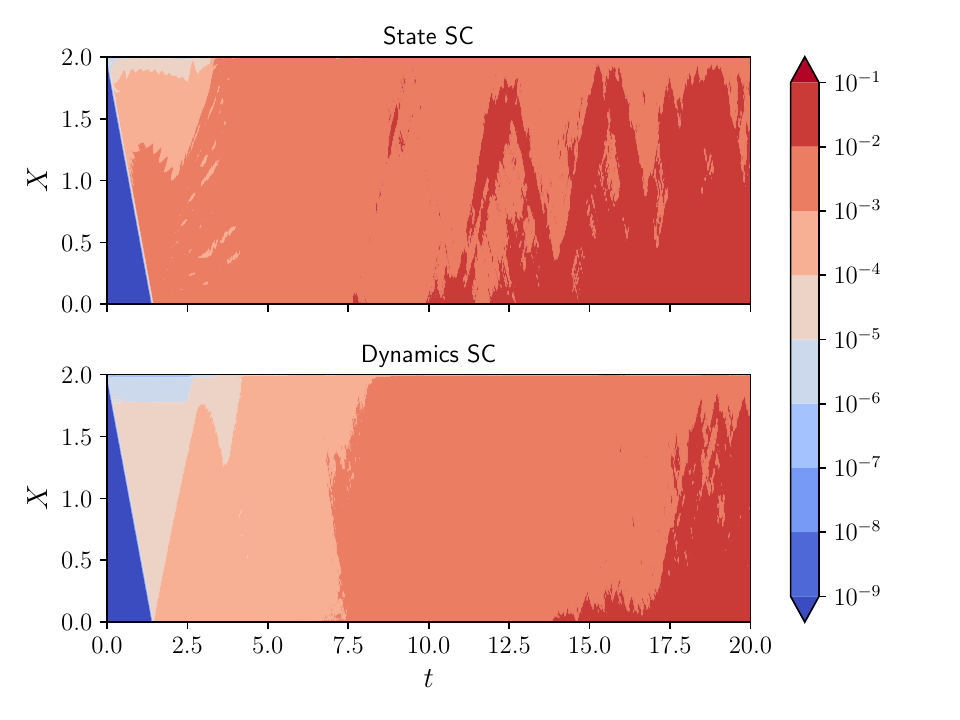}
    \caption{1st order sparse grid, average error}
  \end{subfigure}
  ~
  \begin{subfigure}{0.48\textwidth}
    \centering
    \includegraphics[width=\linewidth]{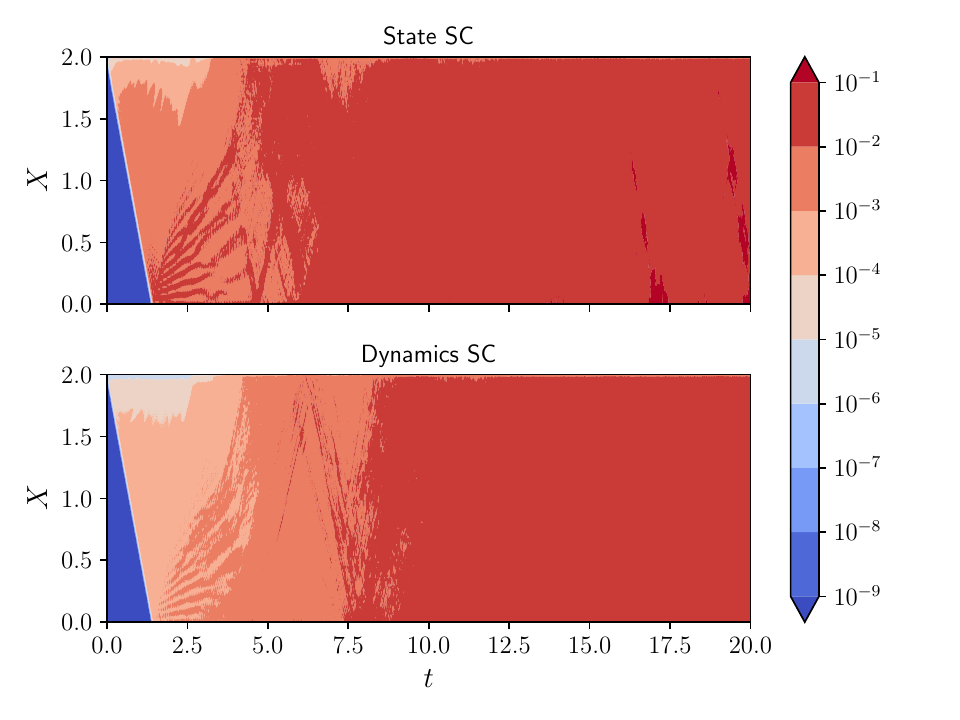}
    \caption{1st order sparse grid, maximum error}
  \end{subfigure}

  \vspace{0.1in}

  \begin{subfigure}{0.48\textwidth}
    \centering
    \includegraphics[width=\linewidth]{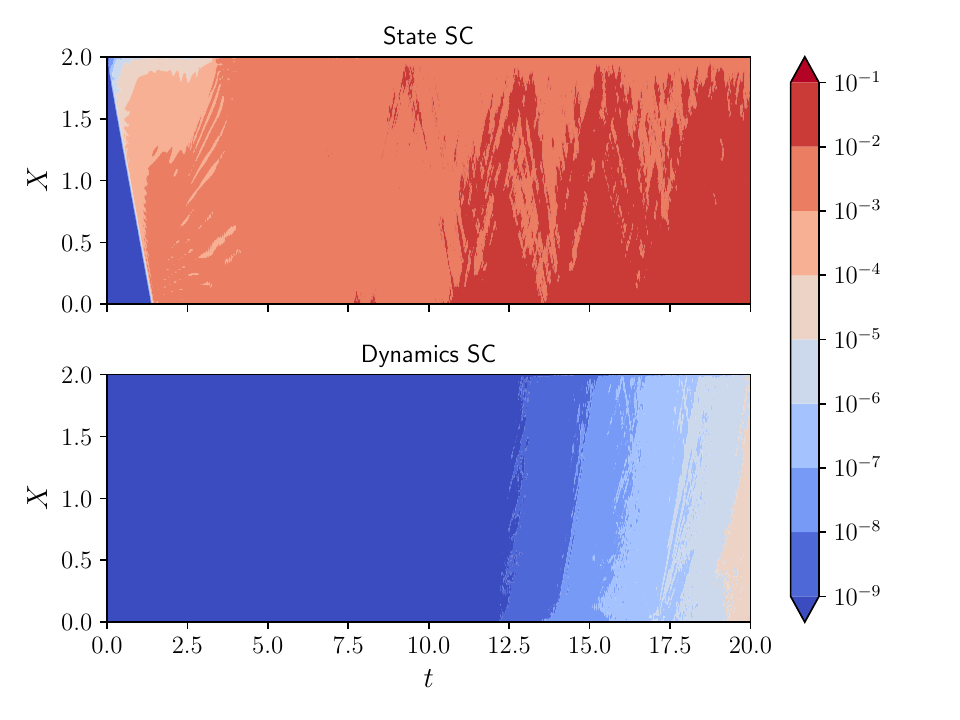}
    \caption{2nd order sparse grid, average error}
  \end{subfigure}
  ~
  \begin{subfigure}{0.48\textwidth}
    \centering
    \includegraphics[width=\linewidth]{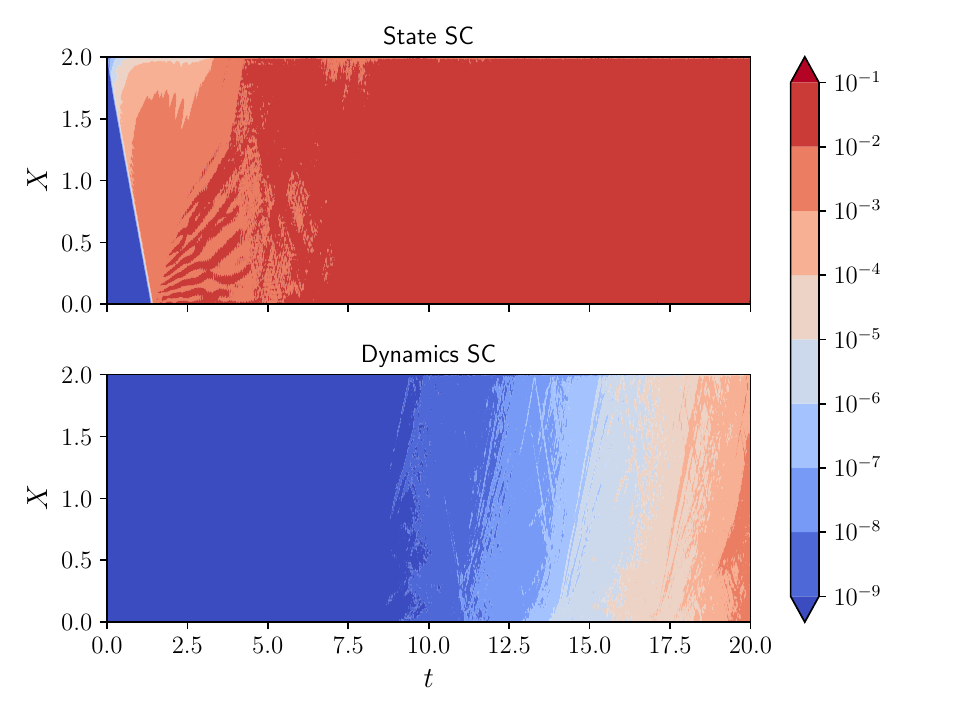}
    \caption{2nd order sparse grid, maximum error}
  \end{subfigure}

  \vspace{0.1in}

  \begin{subfigure}{0.48\textwidth}
    \centering
    \includegraphics[width=\linewidth]{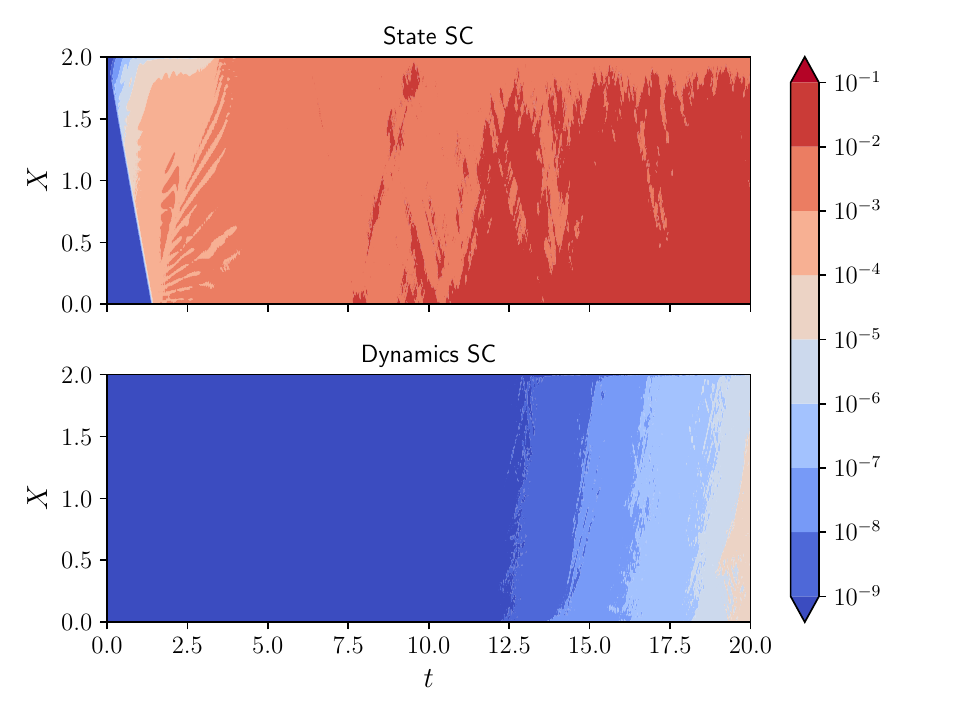}
    \caption{3rd order sparse grid, average error}
  \end{subfigure}
  ~
  \begin{subfigure}{0.48\textwidth}
    \centering
    \includegraphics[width=\linewidth]{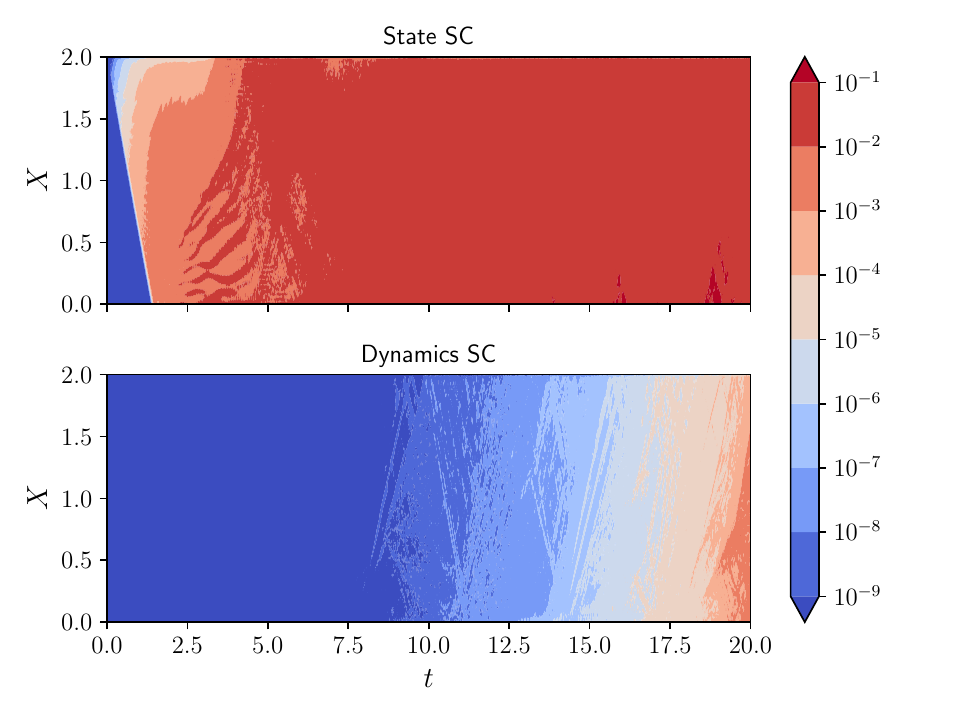}
    \caption{3rd order sparse grid, maximum error}
  \end{subfigure}
  
  \caption{
    Errors in the approximated displacement $u(t, X)$ field using State and Dynamics SC surrogates.
    We plot the point-wise average (left) and maximum (right) errors over 128 sample trajectories using a logarithmic color scale.
    The rows, top to bottom, show the results for 1st, 2nd, and 3rd order sparse grids, corresponding to 5, 13, and 29 exact simulations.
  }
  \label{fig:bar_impact_error_u}
\end{figure*}


%% file: figures/bar_impact/error_v.tex

\begin{figure}[t]
  \centering

  \begin{subfigure}{0.48\textwidth}
    \centering
    \includegraphics[width=\linewidth]{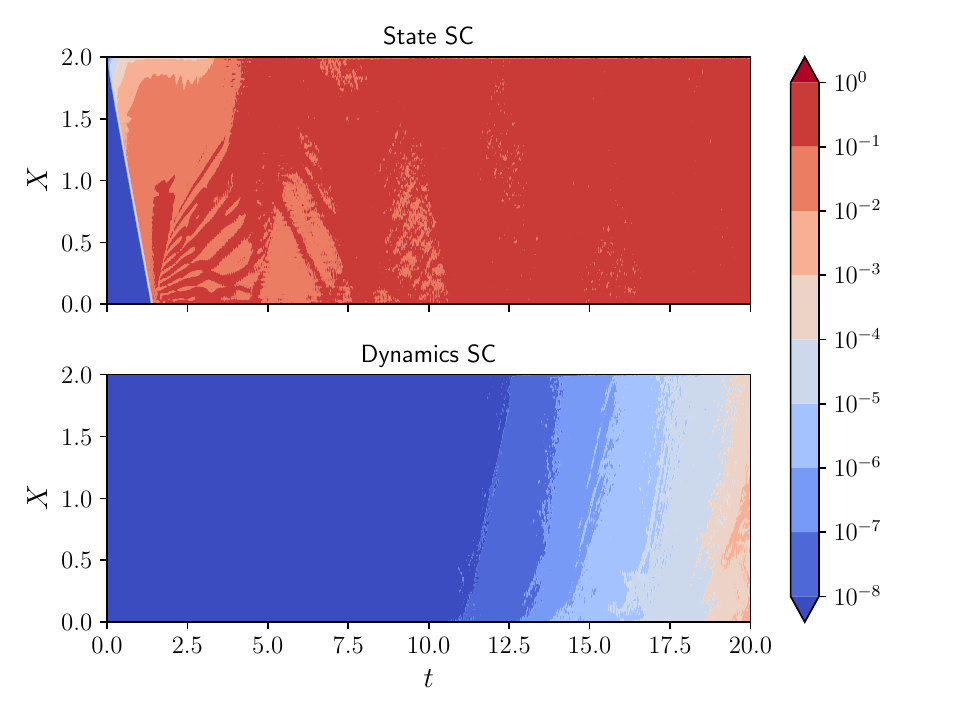}
    \caption{Average error}
  \end{subfigure}
  ~
  \begin{subfigure}{0.48\textwidth}
    \centering
    \includegraphics[width=\linewidth]{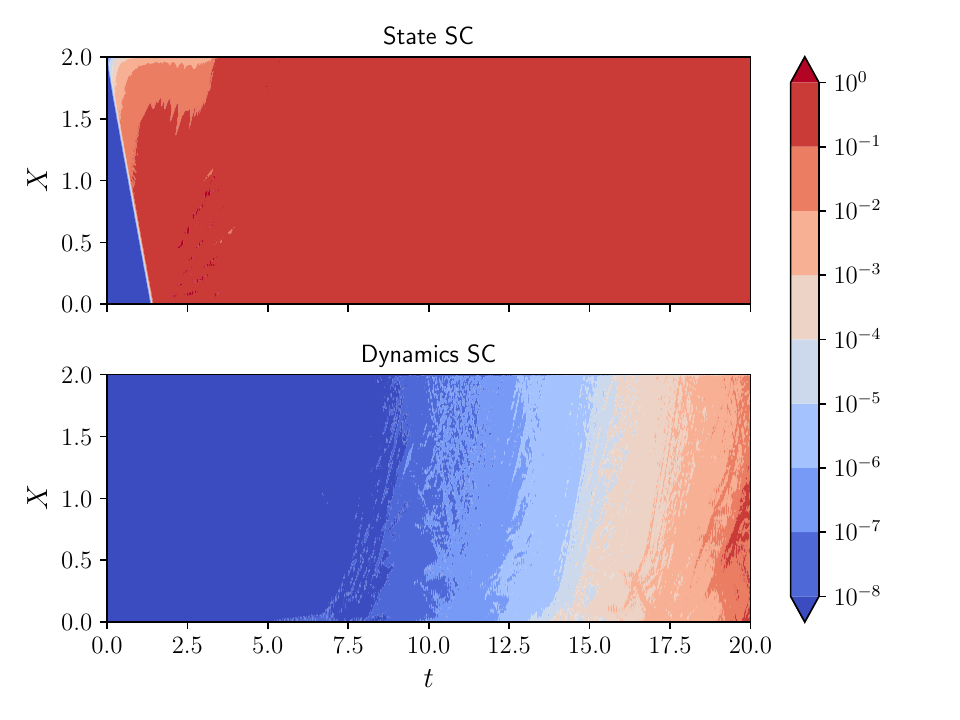}
    \caption{Maximum error}
  \end{subfigure}

  \caption{
    Errors in the approximated velocity fields $\dot{u}(t, X)$ using State and Dynamics SC surrogates with a 2nd order sparse grid.
    We plot the point-wise average (left) and maximum (right) errors over 128 sample trajectories.
  }
  \label{fig:bar_impact_error_v}
\end{figure}


%% file: figures/bar_impact/error_emd.tex

\begin{figure}[t]
  \centering
  \begin{subfigure}{0.48\textwidth}
    \centering
    \includegraphics[width=\linewidth]{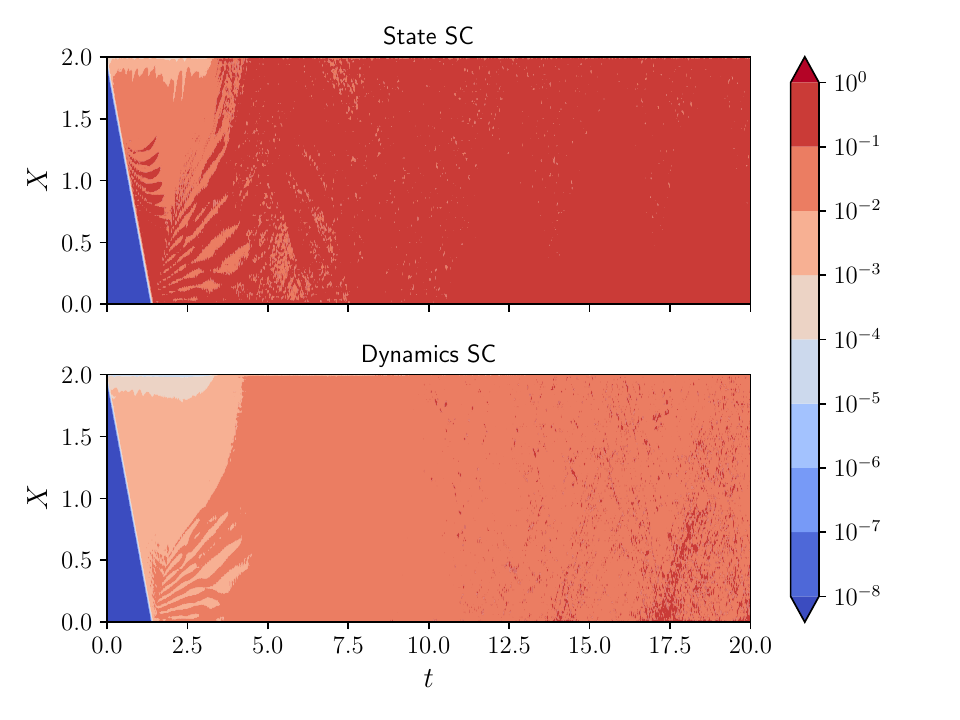}
    \caption{1st order sparse grid}
  \end{subfigure}
  ~
  \begin{subfigure}{0.48\textwidth}
    \centering
    \includegraphics[width=\linewidth]{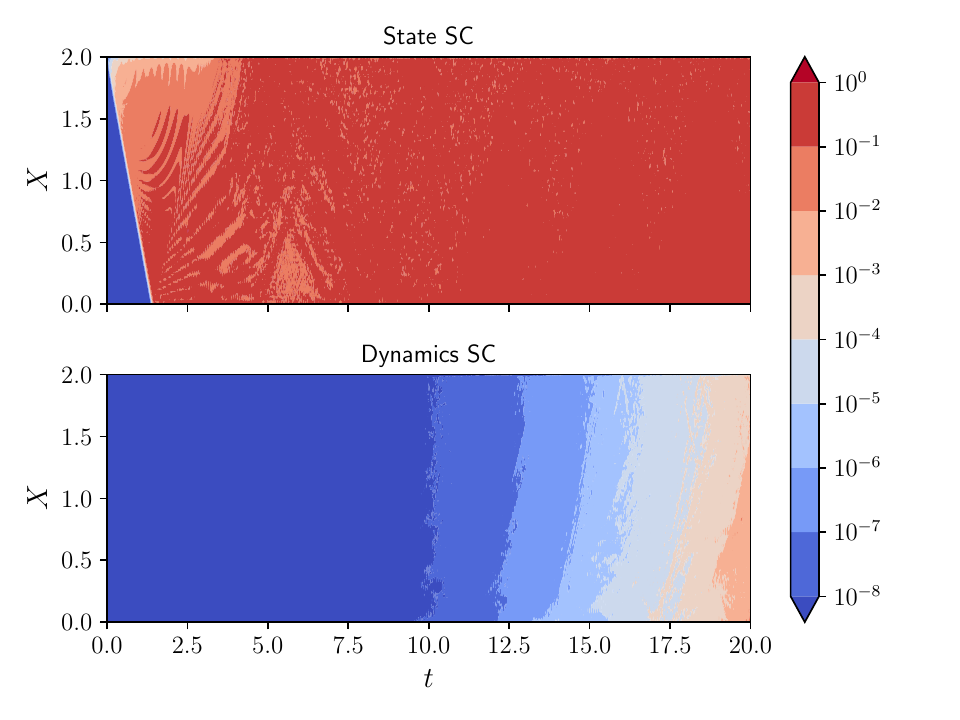}
    \caption{2nd order sparse grid}
  \end{subfigure}
  \caption{
    Wasserstein distance of the empirical state distributions, as estimated from the two SC methods, from the references, which are constructed by running exact simulations.
    The distances are independently computed at each time step and at each spatial node.
  }
  \label{fig:bar_impact_error_emd}
\end{figure}


%% file: figures/notched_plate/setup.tex

\begin{figure}
  \centering
  \includegraphics[width=0.7\textwidth]{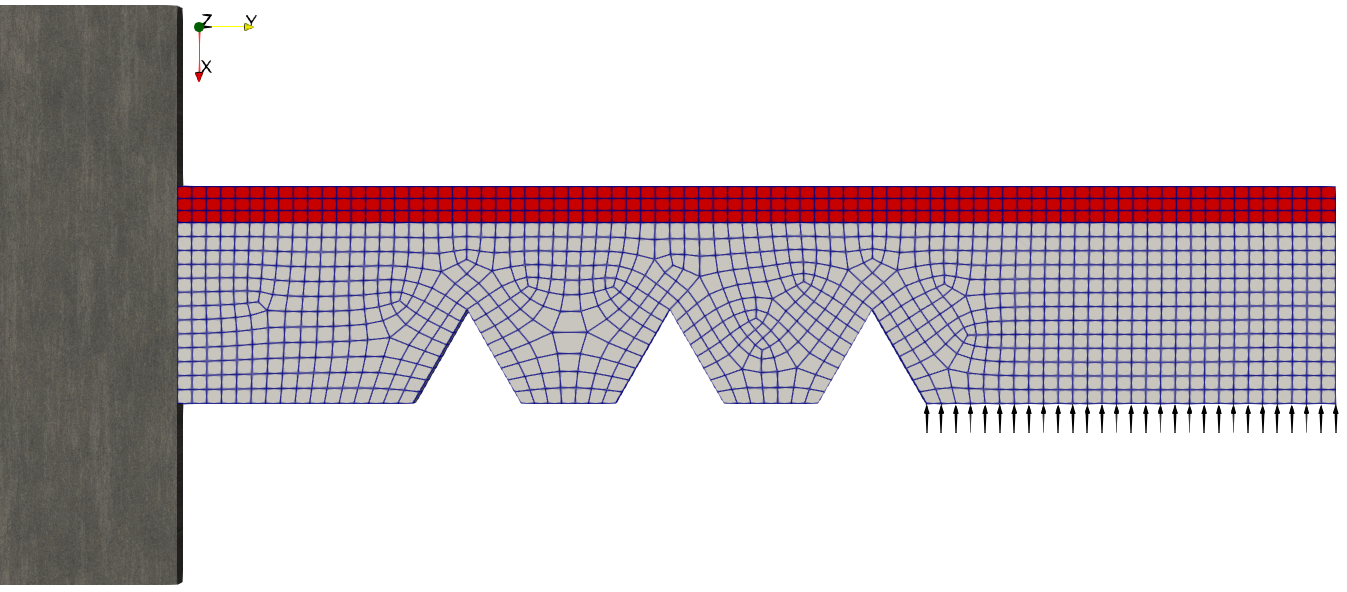}
  \caption{
    Setup for the pseudo-3D notch impact problem.
    The composite plate consists of two materials: glass (in red) and aluminum (in light grey), and the aluminum layer has three notches.
    Left end of the plate is fixed, while the right side of bottom edge is subject to an impact loading, which imparts an initial velocity at the corresponding nodes.
    The structure is discretized using 1169 hexahedral elements, and simulated using open-source Lagrangian solid mechanics software \texttt{NimbleSM}.
  }
  \label{fig:notched_plate_setup}
\end{figure}


%% file: figures/notched_plate/error_u.tex

\begin{figure}[t]
  \centering
  
  \begin{subfigure}{0.48\textwidth}
    \centering
    \includegraphics[width=\linewidth]{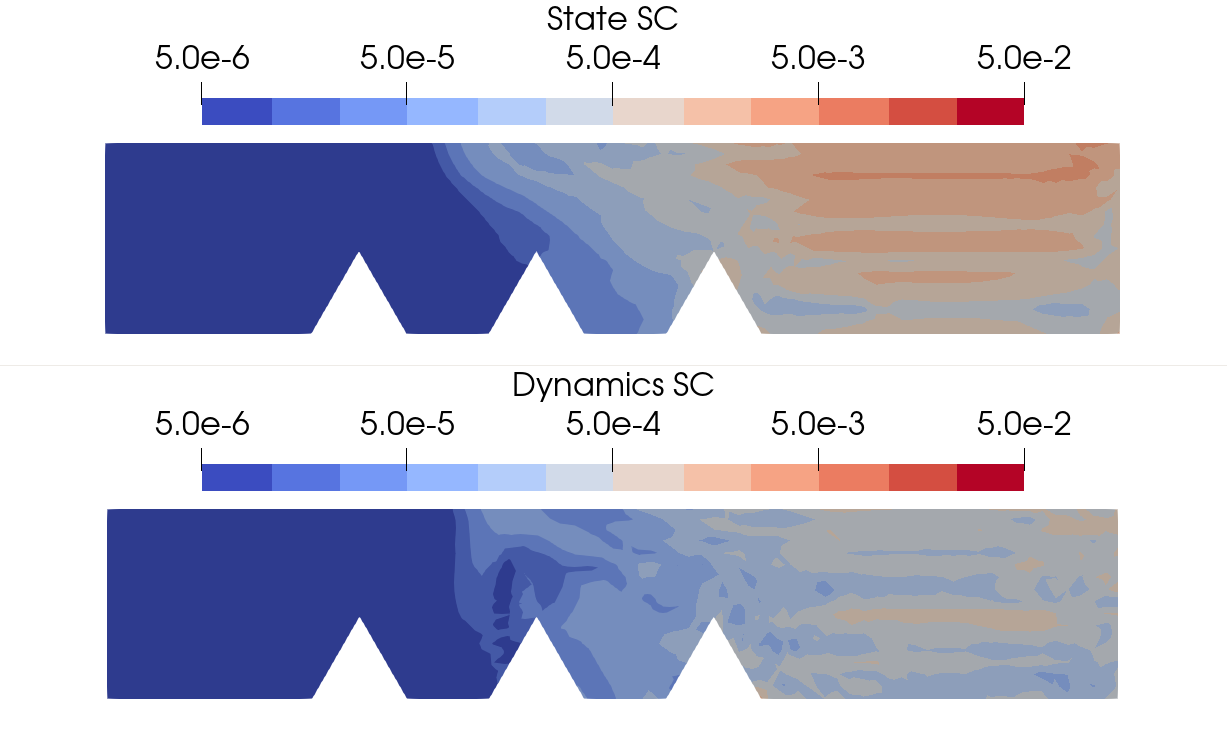}
    \caption{$t = 2.5~\mu$s, average error}
  \end{subfigure}
  ~
  \begin{subfigure}{0.48\textwidth}
    \centering
    \includegraphics[width=\linewidth]{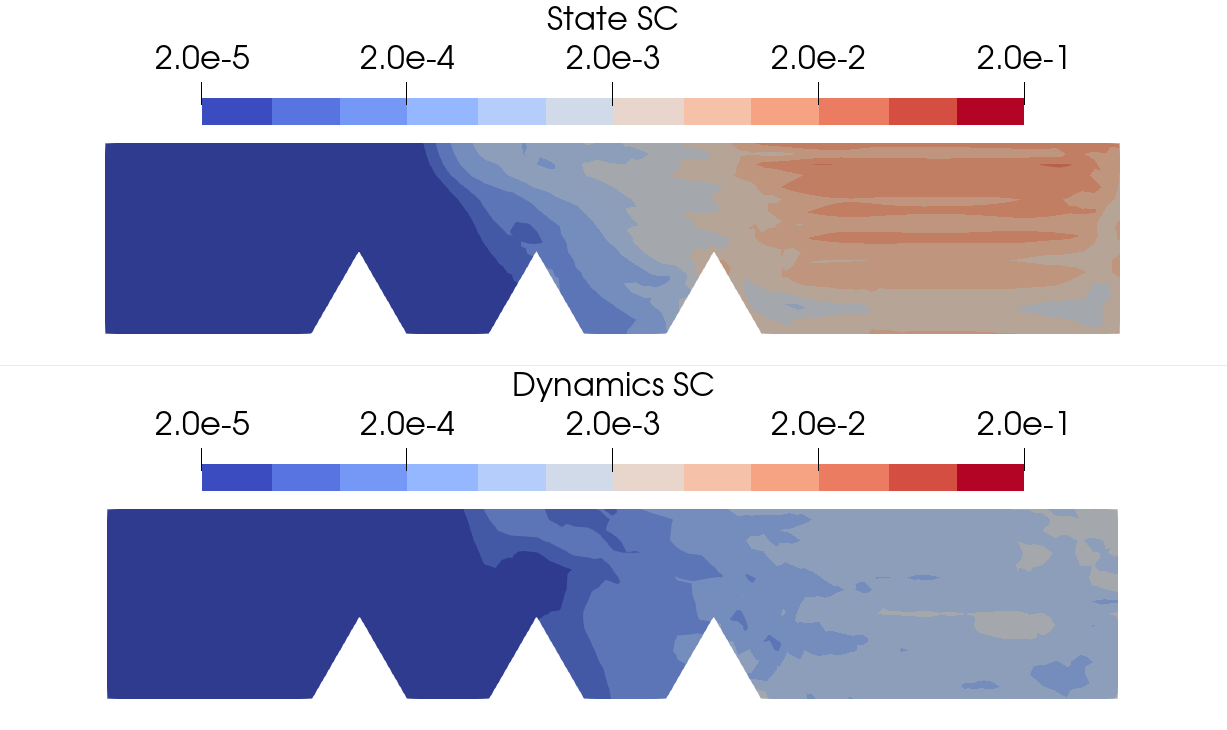}
    \caption{$t = 2.5~\mu$s, maximum error}
  \end{subfigure}

  \bigskip
  
  \begin{subfigure}{0.48\textwidth}
    \centering
    \includegraphics[width=\linewidth]{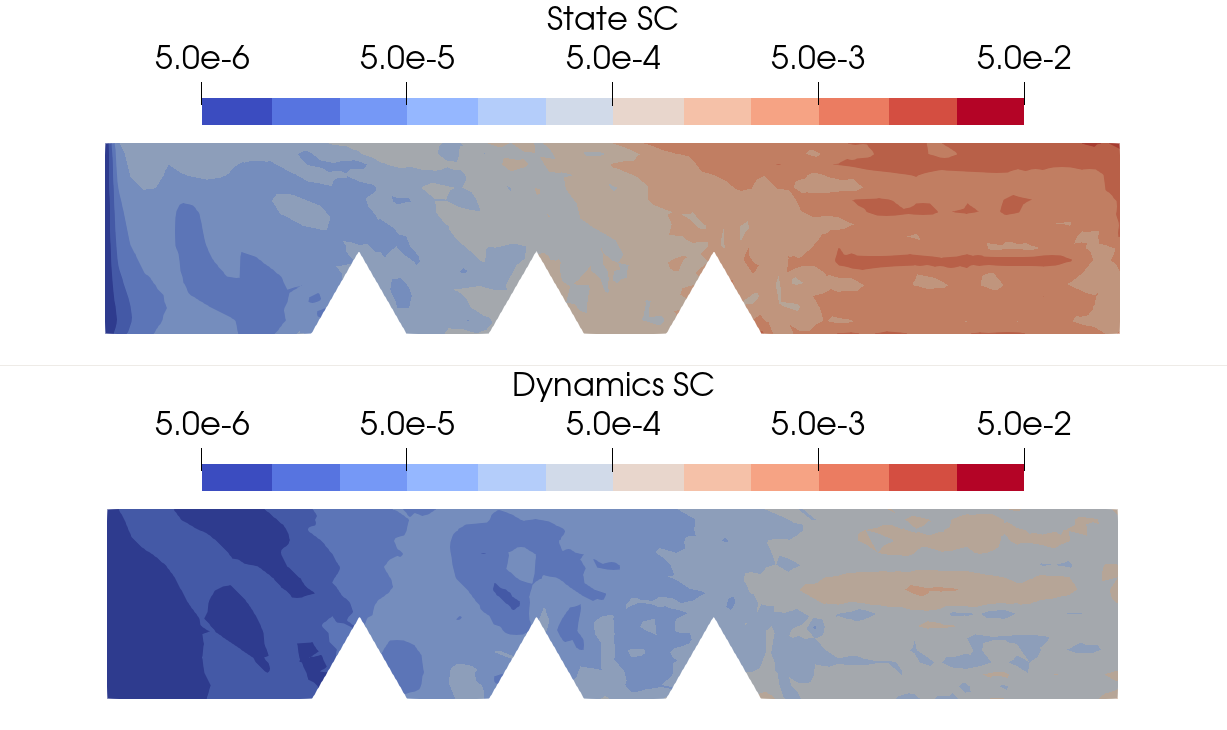}
    \caption{$t = 5.0~\mu$s, average error}
  \end{subfigure}
  ~
  \begin{subfigure}{0.48\textwidth}
    \centering
    \includegraphics[width=\linewidth]{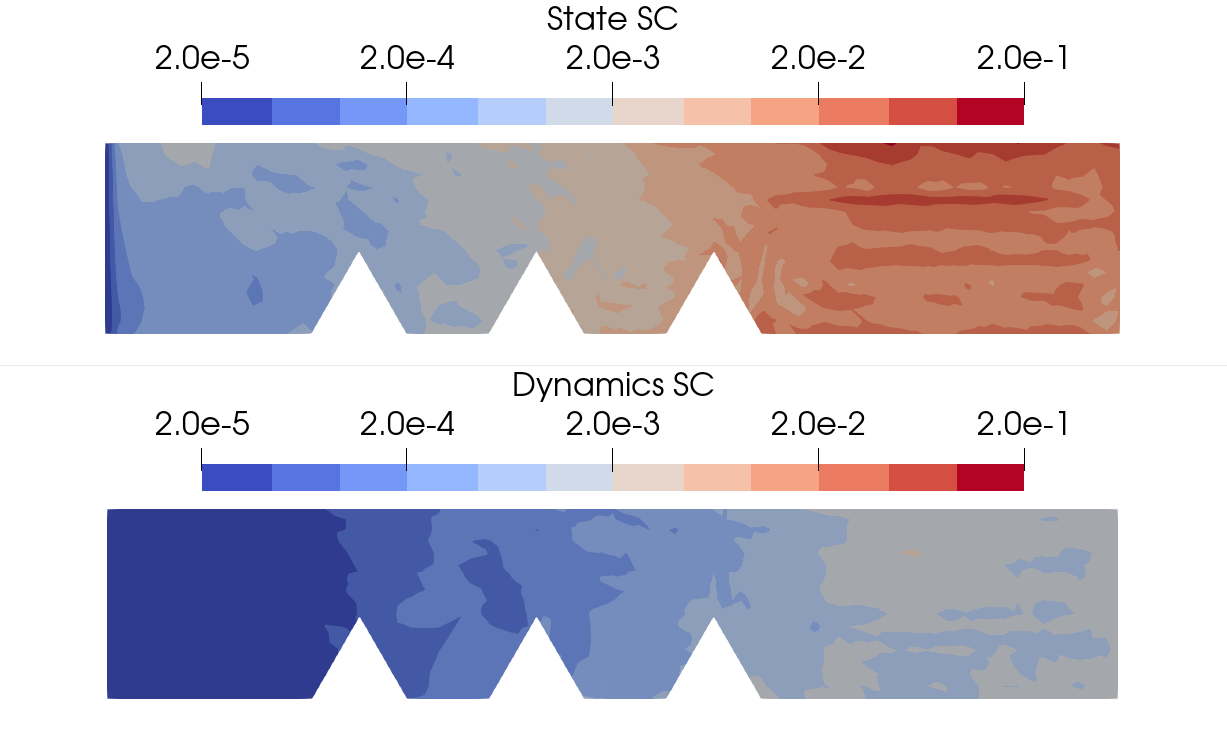}
    \caption{$t = 5.0~\mu$s, maximum error}
  \end{subfigure}

  \bigskip
  
  \begin{subfigure}{0.48\textwidth}
    \centering
    \includegraphics[width=\linewidth]{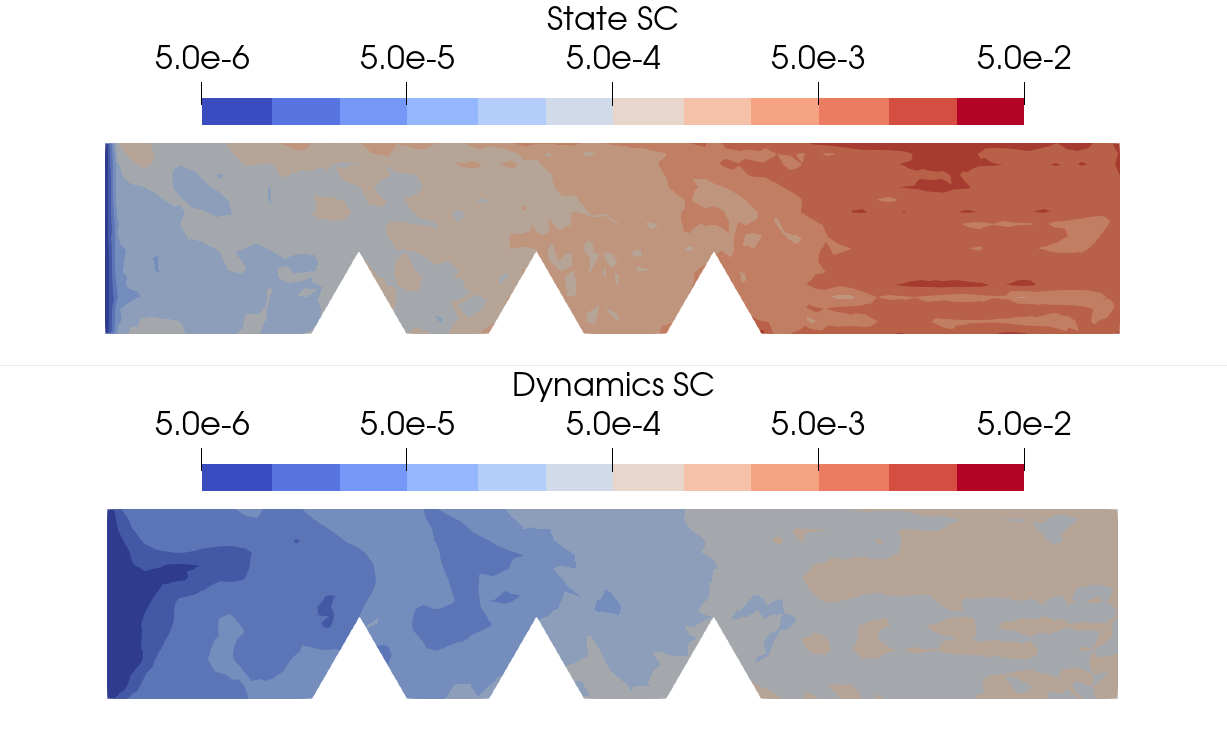}
    \caption{$t = 7.5~\mu$s, average error}
  \end{subfigure}
  ~
  \begin{subfigure}{0.48\textwidth}
    \centering
    \includegraphics[width=\linewidth]{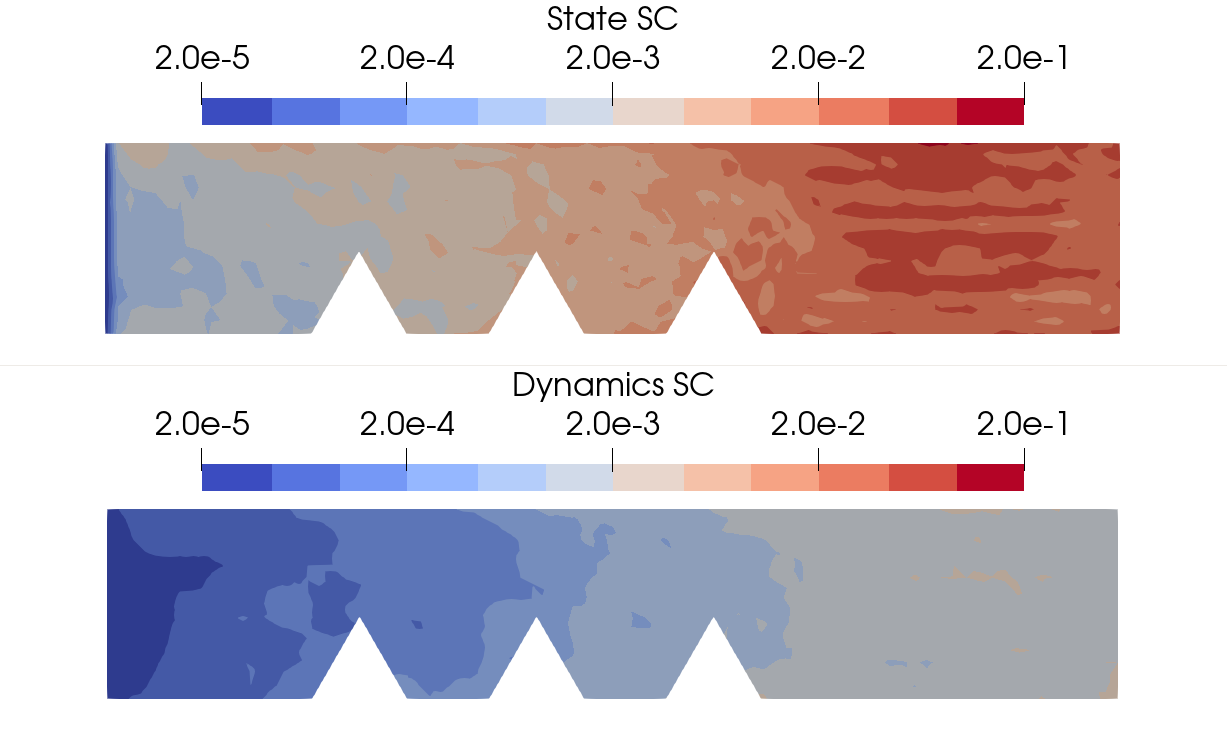}
    \caption{$t = 7.5~\mu$s, maximum error}
  \end{subfigure}

  \caption{
    Comparison of errors in the displacement field $\Vec{u}(t, \Vec{X})$ in the notch impact problem constructed using State and Dynamics SC surrogates.
    The left column shows the error averaged over the 512 sample trajectories, and the right column shows the corresponding maximum error.
  }
  \label{fig:notched_plate_error_u}
\end{figure}


%% file: figures/notched_plate/error_v.tex

\begin{figure}[t]
  \centering
  
  \begin{subfigure}{0.48\textwidth}
    \centering
    \includegraphics[width=\linewidth]{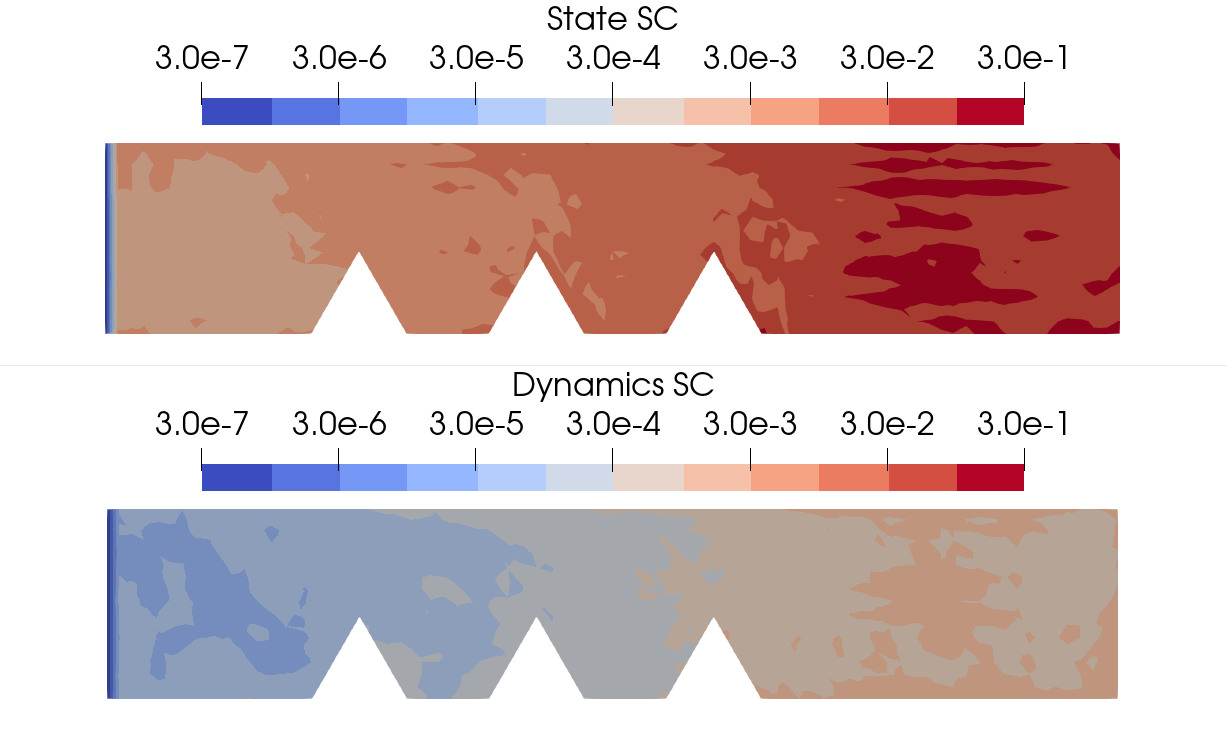}
    \caption{Average error}
  \end{subfigure}
  ~
  \begin{subfigure}{0.48\textwidth}
    \centering
    \includegraphics[width=\linewidth]{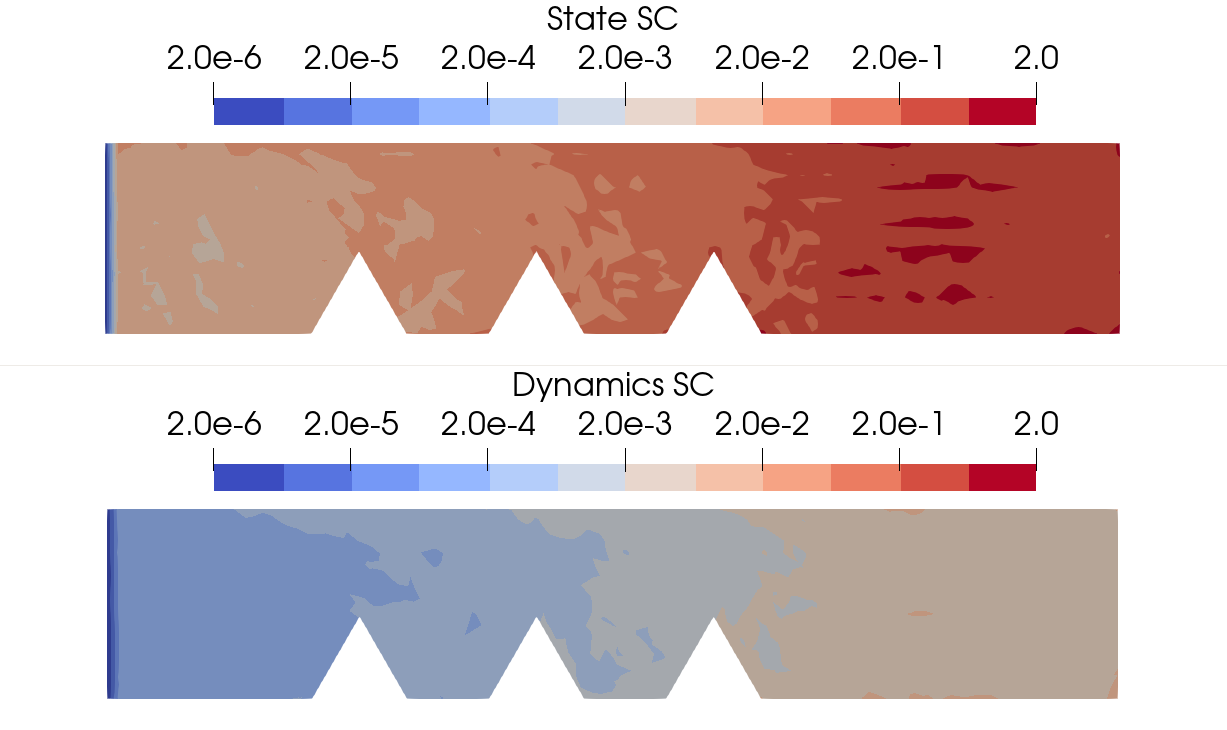}
    \caption{Maximum error}
  \end{subfigure}

  \caption{
    Comparison of errors in the velocity field $\dot{\Vec{u}}(t, \Vec{X})$ at $t = 10~\mu$s in the notch impact problem constructed using State and Dynamics SC surrogates.
    The left panel shows the error averaged over the 512 sample trajectories, and the right panel shows the corresponding maximum error.
  }
  \label{fig:notched_plate_error_v}
\end{figure}


%% file: figures/notched_plate/error_emd.tex

\begin{figure}
  \centering
  
  \begin{subfigure}{0.48\textwidth}
    \centering
    \includegraphics[width=\linewidth]{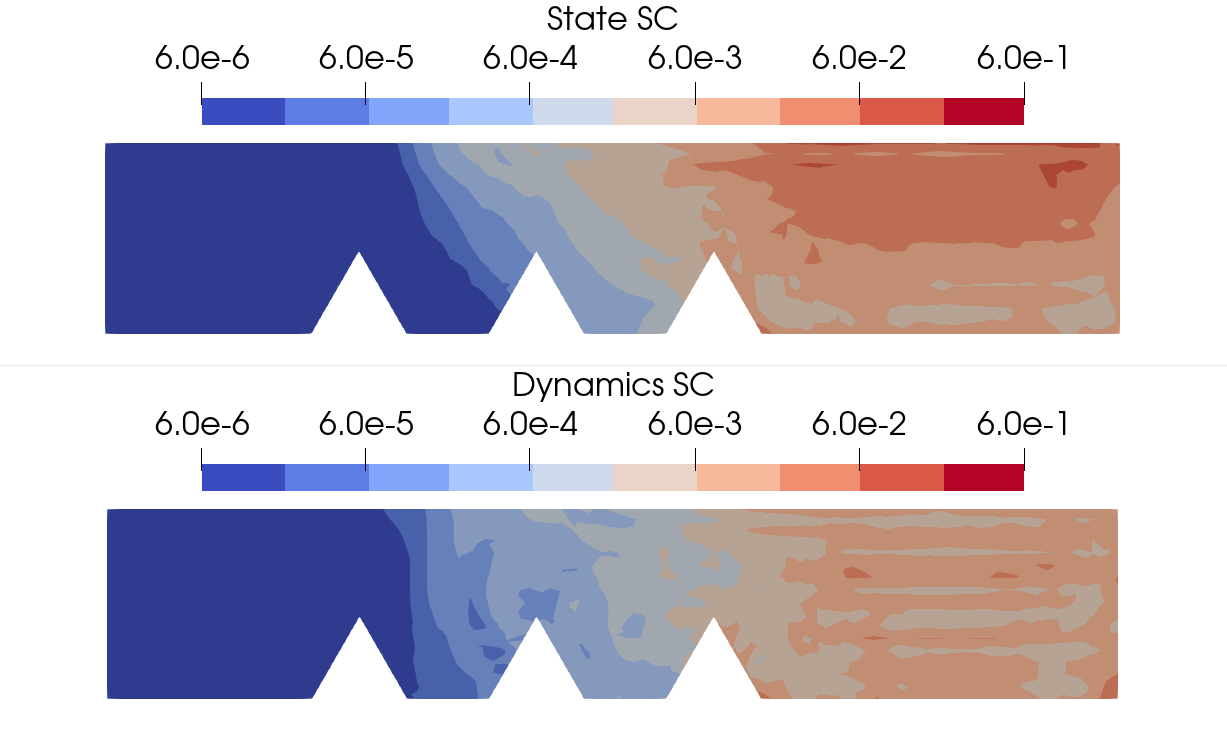}
    \caption{$t = 2.5~\mu$s}
  \end{subfigure}
  ~
  \begin{subfigure}{0.48\textwidth}
    \centering
    \includegraphics[width=\linewidth]{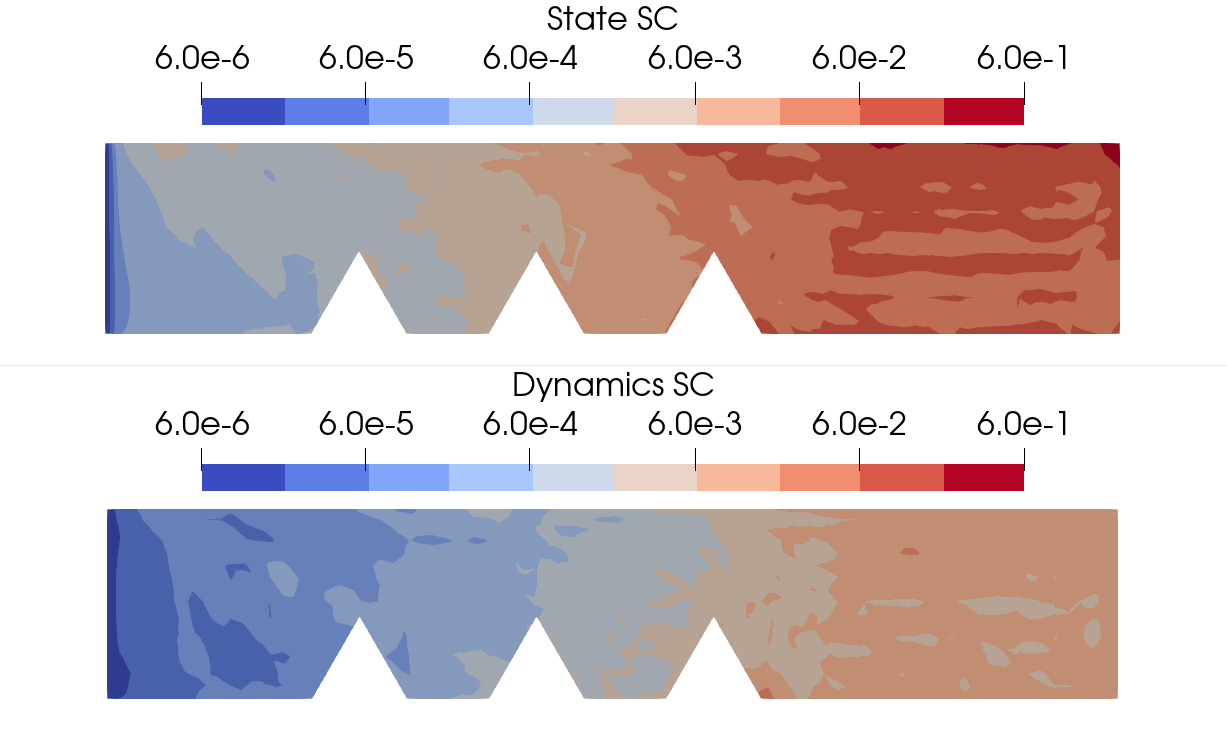}
    \caption{$t = 5.0~\mu$s}
  \end{subfigure}

  \medskip
  
  \begin{subfigure}{0.48\textwidth}
    \centering
    \includegraphics[width=\linewidth]{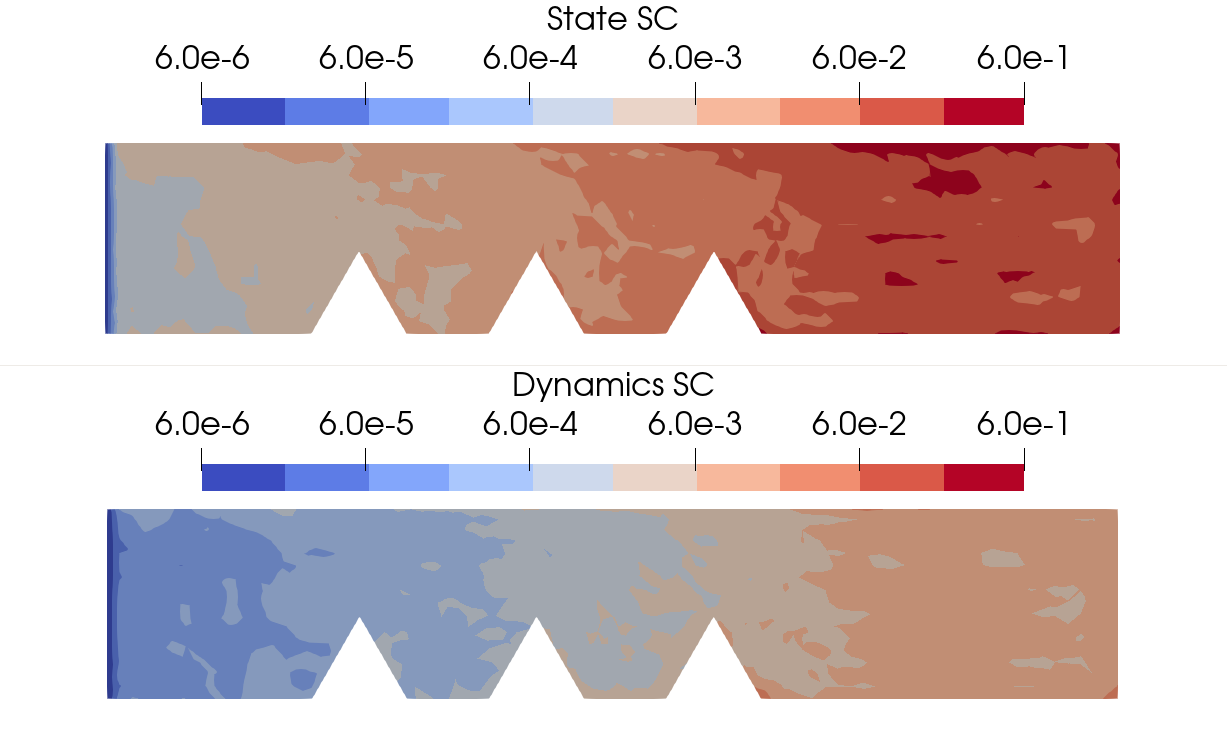}
    \caption{$t = 7.5~\mu$s}
  \end{subfigure}
  ~
  \begin{subfigure}{0.48\textwidth}
    \centering
    \includegraphics[width=\linewidth]{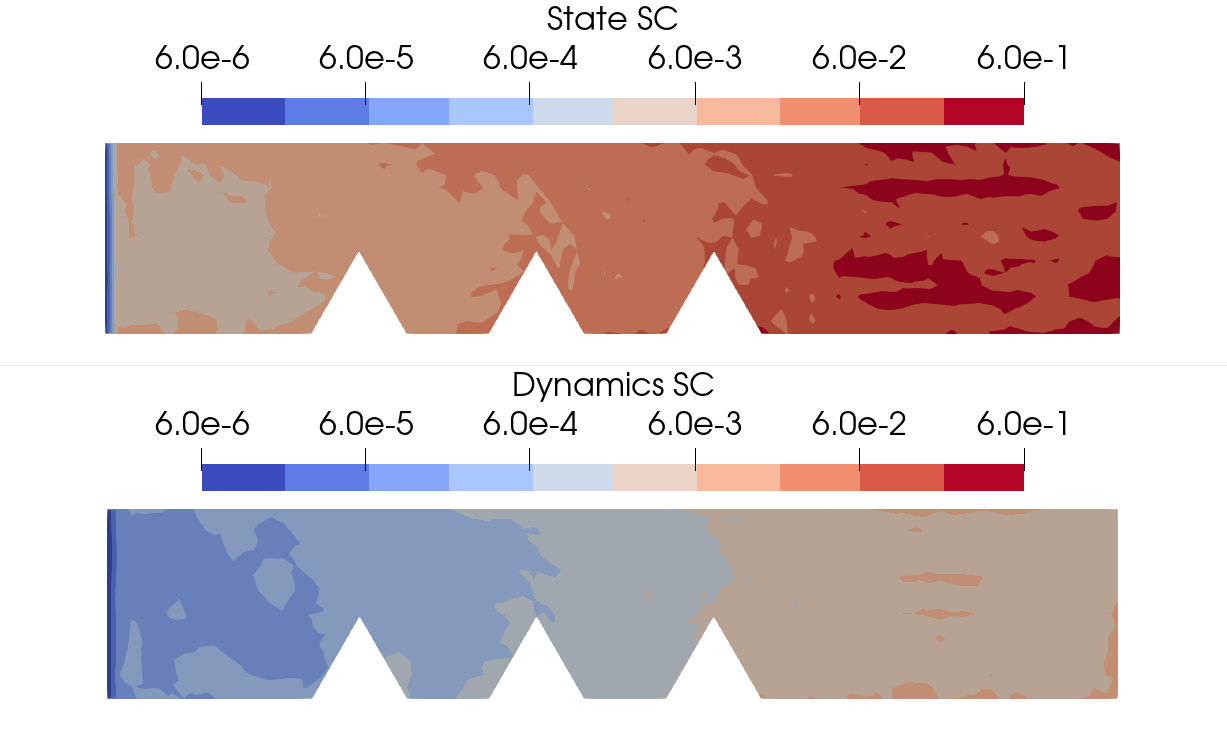}
    \caption{$t = 10.0~\mu$s}
  \end{subfigure}

  \caption{
    Comparison of Wasserstein distance between the empirical state distributions, constructed from the exact simulations and SC surrogate approximations, in the notch impact problem.
    The distances are computed at each time step $t_n$ and at each node $\Vec{X}_j$ of the FEM mesh using the system states $\{(s_u \Vec{u}(t_n, \Vec{X}_j, \Vec{\Lambda}^s_i), s_{\dot{u}} \dot{\Vec{u}}(t_n, \Vec{X}_j, \Vec{\Lambda}^s_i)) : 1 \leq i \leq S\}$ corresponding to the parameter samples $\Vec{\Lambda}^s_i$.
    The factors $s_u = 10^3$ and $s_{\dot{u}} = 10^{-3}$ scale the displacement and velocity fields to the same orders of magnitude; the values are chosen empirically based off the exact simulations.
  }
  \label{fig:notched_plate_error_emd}
\end{figure}


%% file: sections/summary.tex

\section{Summary and Future Work} \label{sec:summary}

UQ is an integral component in the analysis of complex and high-consequence computational codes.
With Monte Carlo ensembles being the method of choice for uncertainty propagation, the need for efficient and accurate reduction of the full model is clear.
Due to its non-intrusive nature and relatively low computational cost, SC is an attractive framework for constructing surrogates; however its accuracy is often limited by the parametric nonlinearity of the QoIs.

We presented an alternate approach, which directly approximates the driving term of a dynamical system.
Our approach is based on the observation that a small number of exactly simulated model trajectories is usually sufficient for learning the form of the model dynamics, as a function of the state vector, over the domain of the uncertain model parameters.
We propose a novel learning framework that combines:
(a) SINDy, which discovers the dependence of the dynamics on the model state, and
(b) SC, which interpolates the dynamics over the stochastic parameter space.
In this way, we create a minimally-invasive approach for uncertainty propagation, requiring only the ability to:
(a) retrieve dynamics evaluations from the target application, and
(b) substitute a learned form of the dynamics into the application code to simulate approximate trajectories.
Consequently we realize order-of-magnitude improvement in the accuracy of the approximate trajectories relative to the widely-used approach of applying SC directly to the state for our test problems.
Furthermore, using the Wasserstein metric to compare empirical state distributions, we demonstrated that our proposed Dynamics SC approach approximates state distributions better than State SC.
We remark that these improvements were expected for the Lorenz system based on the findings in the original SINDy article \cite{brunton2016discovering} and subsequent works.
We demonstrated the same holds for discretized PDE systems, especially for the \texttt{NimbleSM} based notch problem where we approximated a rational stress model with locally linear surrogate.

Dynamics SC does introduce a more complex learning model compared to State SC.
In addition to the application code intrusion outlined in \sref{sec:dynamics_sc}, we also have the added cost of:
(a) managing the history of the exact trajectories at the collocation nodes, and
(b) learning the sparse dynamics coefficients using SINDy.
Learning localized dynamics surrogates from limited trajectory segments within the adaptive semi-online training setup offsets some of this computational complexity.
The choice of SINDy basis functions is also critical: if the dynamics of the system under consideration does not belong to the linear span of the selected basis functions, then the approximation can become very poor.

In future work, we will address a limitation  of using SINDy as the dynamics discovery framework:
namely, the number of candidate functions used in dynamics approximation grow rapidly for high-dimensional problems.
We propose to use functional tensor formats, such as the function train \cite{gorodetsky2018gradientbased, gorodetsky2019continuous}, to break this curse of dimensionality.
We also wish to tackle more complex material behavior in solid mechanics, such as plasticity, with the Dynamics SC approach.
